# RATES OF CONVERGENCE OF A TRANSIENT DIFFUSION IN A SPECTRALLY NEGATIVE LÉVY POTENTIAL

### By Arvind Singh

### *Université Pierre et Marie Curie*

We consider a diffusion process $X$ in a random Lévy potential $\mathbb{V}$ which is a solution of the informal stochastic differential equation

$$\begin{cases} dX_t = d\beta_t - \frac{1}{2}\mathbb{V}'(X_t)\,dt, \\ X_0 = 0, \end{cases}$$

($\beta$ B. M. independent of $\mathbb{V}$). We study the rate of convergence when the diffusion is transient under the assumption that the Lévy process $\mathbb{V}$ does not possess positive jumps. We generalize the previous results of Hu–Shi–Yor for drifted Brownian potentials. In particular, we prove a conjecture of Carmona: provided that there exists $0 < \kappa < 1$ such that $\mathbf{E}[e^{\kappa \mathbb{V}_1}] = 1$, then $X_t/t^\kappa$ converges to some nondegenerate distribution. These results are in a way analogous to those obtained by Kesten–Kozlov–Spitzer for the transient random walk in a random environment.

**1. Introduction.** Let $(\mathbb{V}(x), x \in \mathbb{R})$ be a càdlàg, real-valued stochastic process with $\mathbb{V}(0) = 0$, defined on some probability space $(\Omega, \mathbf{P})$. We consider a diffusion process $X$, solution of the informal stochastic differential equation

$$\begin{cases} dX_t = d\beta_t - \frac{1}{2}\mathbb{V}'(X_t)\,dt, \\ X_0 = 0, \end{cases}$$

where $(\beta_s, s \geq 0)$ is a standard Brownian motion independent of $\mathbb{V}$. Formally, one can see $X$ as a diffusion process whose conditional generator, given $\mathbb{V}$, is

$$\frac{1}{2}e^{\mathbb{V}(x)}\frac{d}{dx}\left(e^{-\mathbb{V}(x)}\frac{d}{dx}\right).$$

We call $X$ a diffusion in the random potential $\mathbb{V}$. Somehow, this process may be thought as the continuous analogue of the random walk in random

---

environment (see Schumacher [18] or Shi [20] for a connection between the two models). In particular, both models exhibit similar interesting features such as asymptotic sub-linear growth.

For instance, if $\mathbb{V}$ is a two-sided Brownian motion, then $X$ is recurrent and Brox [4] proved an equivalent of Sinai's theorem [21] for random walk in a random environment, that is, $X_t / \log^2 t$ converges to some nondegenerate distribution as $t$ goes to infinity.

When the potential process $\mathbb{V}$ is a drifted Brownian motion ($\mathbb{V}_x = B_x - \frac{\kappa}{2}x$, with $\kappa > 0$ and $B$ a two-sided Brownai motion), the diffusion is transient toward $+\infty$. More precisely, Kawazu and Tanaka [12] showed that the rate of convergence to infinity depends on the value of $\kappa$:

- If $0 < \kappa < 1$, then $\frac{1}{t^\kappa} X_t$ converges in law, as $t$ goes to infinity, toward a nondegenerate positive random variable.
- If $\kappa = 1$, then $\frac{\log t}{t} X_t$ converges in probability toward $\frac{1}{4}$.
- If $\kappa > 1$, then $\frac{1}{t} X_t$ converges almost surely toward $\frac{\kappa - 1}{4}$.

Refined results were later obtained by Tanaka [22] and Hu, Shi and Yor [11], in particular, they proved a central-limit type theorem when $\kappa > 1$. We point out that these results are the analogue, when the potential is a drifted Brownian motion, of those previously obtained by Kesten, Kozlov and Spitzer [14] for the discrete model of the random walk in a random environment. However, the results of Kesten, Kozlov and Spitzer hold for a wide class of environments, whereas few results are available in the continuous setting for general potentials. One would certainly like to extend the results of [11] and [22] for drifted Brownian motion to a wider class of potentials. In this spirit, Carmona [5] considered the case where $\mathbb{V}$ is a two-sided Lévy process and proved, by use of ergodic theorems that, if $\Phi$ denotes the Laplace exponent of $\mathbb{V}$,

$$(1.1) \qquad \mathbf{E}[e^{\lambda \mathbb{V}_t}] = e^{t\Phi(\lambda)}, \qquad t \geq 0, \lambda \in \mathbb{R}$$

[note that $\Phi(\lambda)$ may be infinite], then:

- If $\Phi(1) < 0$, then $X_t / t$ converges almost surely, as $t$ goes to infinity, toward $-\Phi(1)/2$.
- If $\Phi(-1) < 0$ then $X_t / t$ converges almost surely, as $t$ goes to infinity, toward $\Phi(-1)/2$.
- Otherwise, $X_t / t$ converges almost surely toward 0.

Carmona also conjectured that when the limiting velocity is zero, assuming that there exists $0 < \kappa < 1$ such that $\Phi(\kappa) = 0$, then one should observe the same asymptotic behavior as in the case of a drifted Brownian potential, that is, the rate of growth of $X_t$ should again be of order $t^\kappa$. We prove that this is the case when $\mathbb{V}$ is a spectrally negative Lévy process (i.e., a Lévy process without positive jumps).

Throughout this paper we will make the following assumption on $\mathbb{V}$:

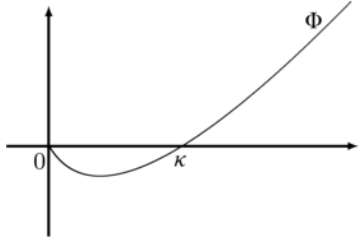

Fig. 1.   *The Laplace exponent $\Phi$.*

Assumption 1.1.   The following hold:

(a)  $(\mathbb{V}_x, x \in \mathbb{R})$ is a càdlàg locally bounded process with $\mathbb{V}_0 = 0$ and the two processes $(\mathbb{V}_x, x \geq 0)$ and $(\mathbb{V}_{-x}, x \geq 0)$ are independent.

(b)  $(\mathbb{V}_x, x \geq 0)$ is a Lévy process with no positive jumps which is not the opposite of a subordinator and is such that $\lim_{x \to \infty} \mathbb{V}_x = -\infty$ almost surely.

(c)  $(\mathbb{V}_{-x}, x \geq 0)$ is such that $\int_0^\infty e^{\mathbb{V}_{-x}} \, dx = \infty$ almost surely.

Let us first make some comments concerning our assumptions.

– Note that (c) is a weak condition. For instance, it is fulfilled whenever $(\mathbb{V}_{-x}, x \geq 0)$ is a Lévy process which does not diverge to $-\infty$. In fact, (c) is only to ensure that the diffusion $X$ does not go to $-\infty$ with positive probability. Otherwise, we are not really concerned about the behavior of $\mathbb{V}$ for negative $x$'s. In particular, the process $(\mathbb{V}_{-x}, x \geq 0)$ may have jumps of both signs. See Figure 2.

– Since $(\mathbb{V}_x, x \geq 0)$ has no positive jumps, its Laplace exponent $\Phi$ given by (1.1) is finite at least for all $\lambda \in [0, \infty)$. The assumption that $\mathbb{V}$ is not the opposite of a subordinator implies that $\Phi(\lambda) \to \infty$ as $\lambda \to \infty$. Moreover, since $\mathbb{V}$ is transient toward $-\infty$, the right derivative of $\Phi$ at $0+$ is such that $\Phi'(0+) = \mathbf{E}[\mathbb{V}_1] \in [-\infty, 0)$. Thus, the strict convexity of $\Phi$ implies that $\mathbb{V}$ fulfills the so-called Cramér's condition (see Figure 1): there exists a unique $\kappa > 0$ such that

(1.2)                              $\Phi(\kappa) = 0.$

In particular, $\Phi(x) < 0$ for all $x \in (0, \kappa)$, whereas $\Phi(x) > 0$ for all $x > \kappa$.

We introduce the scale function of the diffusion $X$:

(1.3)                 $A(x) \overset{\text{def}}{=} \int_0^x e^{\mathbb{V}_y} \, dy \qquad \text{for } x \in [-\infty, \infty].$

On the one hand, Assumption 1.1(c) implies that

(1.4)                 $\lim_{x \to -\infty} A(x) = A(-\infty) = -\infty, \qquad \mathbf{P}\text{-a.s.}$

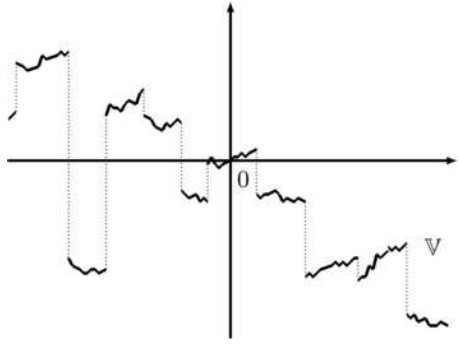

Fig. 2.    *Sample path of $\mathbb{V}$.*

On the other hand, in view of Assumption 1.1(b), for $0 < \delta < -\mathbf{E}[\mathbb{V}_1]$, the Lévy process $(\mathbb{V}_x + \delta x, x \geq 0)$ also diverges toward $-\infty$. This entails

$$(1.5) \qquad \lim_{x \to +\infty} A(x) = A(+\infty) < \infty, \qquad \mathbf{P}\text{-a.s.}$$

Combining (1.4) and (1.5), it is easy to check that $X$ is transient toward $+\infty$ (see [20] for details). We now introduce the hitting time of level $r \geq 0$ for the diffusion:

$$(1.6) \qquad H(r) \overset{\text{def}}{=} \inf\{t \geq 0, X_t = r\}.$$

Let $\mathcal{N}$ stand for a Gaussian $\mathcal{N}(0,1)$ variable. For $\alpha \in (0,1) \cup (1,2)$, let $\mathcal{S}_\alpha^{\text{ca}}$ be a completely asymmetric stable variable with characteristic function

$$\mathbf{E}[e^{it\mathcal{S}_\alpha^{\text{ca}}}] = \exp\left(-|t|^\alpha\left(1 - i\operatorname{sgn}(t)\tan\left(\frac{\pi\alpha}{2}\right)\right)\right)$$

($S_\alpha^{\text{ca}}$ is positive when $\alpha < 1$). Let $\mathcal{C}^{\text{ca}}$ also denote a completely asymmetric Cauchy variable with characteristic function

$$\mathbf{E}[e^{it\mathcal{C}^{\text{ca}}}] = \exp\left(-\left(|t| + it\frac{2}{\pi}\log|t|\right)\right).$$

We can now state our main theorem.

THEOREM 1.1.    *Recall that $\kappa$ defined by (1.2) is the unique positive root of the Laplace exponent $\Phi$ of $\mathbb{V}$. We denote by $\Phi'$ the derivative of $\Phi$. Set*

$$\mathbf{K} \overset{\text{def}}{=} \mathbf{E}\left[\left(\int_0^\infty e^{\mathbb{V}_y}\,dy\right)^{\kappa-1}\right].$$

*This constant (which only depends on the potential $\mathbb{V}$) is finite. When $\kappa > 1$ [i.e., when $\Phi(1) < 0$], set $\mathbf{m} \overset{\text{def}}{=} -2/\Phi(1) > 0$. We have the following, depending on the value of $\kappa$:*

(a) If $0 < \kappa < 1$,

$$\frac{1}{r^{1/\kappa}} H(r) \xrightarrow[r \to \infty]{\text{law}} 2 \left( \frac{\pi \kappa^2 \mathbf{K}^2}{2 \sin(\pi \kappa / 2) \Phi'(\kappa)} \right)^{1/\kappa} \mathcal{S}_\kappa^{\text{ca}}.$$

(b) If $\kappa = 1$, there exists a function $f$ with $f(r) \sim \frac{2}{\Phi'(1)} r \log r$ such that

$$\frac{1}{r} (H(r) - f(r)) \xrightarrow[r \to \infty]{\text{law}} \left( \frac{\pi}{\Phi'(1)} \right) \mathcal{C}^{\text{ca}}.$$

(c) If $1 < \kappa < 2$,

$$\frac{1}{r^{1/\kappa}} (H(r) - \mathbf{m} r) \xrightarrow[r \to \infty]{\text{law}} 2 \left( \frac{\pi \kappa^2 \mathbf{K}^2}{2 \sin(\pi \kappa / 2) \Phi'(\kappa)} \right)^{1/\kappa} \mathcal{S}_\kappa^{\text{ca}}.$$

(d) If $\kappa = 2$,

$$\frac{1}{\sqrt{r \log r}} (H(r) - \mathbf{m} r) \xrightarrow[r \to \infty]{\text{law}} \left( \frac{-4}{\Phi(1) \sqrt{\Phi'(2)}} \right) \mathcal{N}.$$

(e) If $\kappa > 2$,

$$\frac{1}{\sqrt{r}} (H(r) - \mathbf{m} r) \xrightarrow[r \to \infty]{\text{law}} \sqrt{\frac{8 (\Phi(2) - 4 \Phi(1))}{\Phi(1)^3 \Phi(2)}} \mathcal{N}.$$

This theorem gives precise asymptotics for $H(r)$. It is well known that these estimates may in turn be used to obtain asymptotics for $X_t$, $\sup_{s \le t} X_s$ and $\inf_{s \ge t} X_s$ (see [12] for details). For example, when $0 < \kappa < 1$, (a) of the theorem entails

$$\frac{X_t}{t^\kappa} \xrightarrow[t \to \infty]{\text{law}} \frac{2^{1-\kappa} \sin(\pi \kappa / 2) \Phi'(\kappa)}{\pi \kappa^2 \mathbf{K}^2} \left( \frac{1}{\mathcal{S}_\kappa^{\text{ca}}} \right)^\kappa.$$

The same result also holds for $\sup_{s \le t} X_s$ or $\inf_{s \ge t} X_s$ in place of $X_t$.

One would certainly wish to express the value of the constant $\mathbf{K}$ in terms of the characteristics of the Lévy process $\mathbb{V}$. Although there is to our knowledge no explicit formula for this constant, there are a few cases where the calculations may be carried to their full extent.

EXAMPLE 1.1. We consider a potential of the form $\mathbb{V}_x = B_x - \frac{\kappa}{2} x$ with $\kappa > 0$ and where $B$ is a two-sided standard Brownian motion. According to Dufresne [8] (see also Proposition 2.2 of [7]), the random variable $\int_0^\infty e^{\mathbb{V}_s} ds$ has the same law as $\frac{2}{\gamma_\kappa}$, where $\gamma_\kappa$ denotes a gamma variable with parameter $\kappa$. Therefore, the constant $\mathbf{K}$ may be explicitly calculated:

$$\mathbf{K} = \frac{2^{\kappa-1}}{\Gamma(\kappa)}$$

($\Gamma$ denotes Euler's Gamma function). Thus, we recover the results of Hu, Shi and Yor [11] and Tanaka [22], except for $\kappa = 1$ where we do not have the explicit form of the centering function $f$.

EXAMPLE 1.2. We consider a potential of the form

$$\mathbb{V}_x = cx - \tau_x \qquad \text{for } x \geq 0,$$

with $c > 0$ and where $\tau$ is a subordinator without drift and whose Lévy measure $\nu$ has the form $\nu[x, \infty) = ae^{-bx}$ with $a, b > 0$. Then, the Laplace exponent of $\mathbb{V}$ is given by

$$\Phi(\lambda) = c\lambda - \frac{a\lambda}{\lambda + b} \qquad \text{for all } \lambda \geq 0.$$

Since $\mathbf{E}[\mathbb{V}_1] = c - \frac{a}{b}$, Assumption 1.1 is fulfilled whenever $c < \frac{a}{b}$, in which case Theorem 1.1 holds with $\kappa = \frac{a}{c} - b$. According to Proposition 2.1 of [6], the density $k$ of the integral functional $\int_0^\infty e^{\mathbb{V}_x} dx$ satisfies the differential equation

$$(1 + cx)k(x) = a \int_x^\infty \left(\frac{x}{u}\right)^b k(u) \, du.$$

This equation may be explicitly solved and we find

$$k(x) = \left(\frac{c^{b+1}\Gamma(a/c+1)}{\Gamma(a/c-b)\Gamma(b+1)}\right) \frac{x^b}{(1+cx)^{1+a/c}}.$$

Thus, we can again calculate the value of the constant of Theorem 1.1,

$$\mathbf{K} = \int_0^\infty x^{\kappa-1} k(x) \, dx = \frac{\Gamma(a/c)}{c^{a/c-b-1}\Gamma(a/c-b)\Gamma(b+1)}.$$

In the case of a drifted Brownian potential, in order to obtain the rates of transience of the diffusion, Kawazu and Tanaka [12] and Tanaka [22] made use of Kotani's formula, whereas Hu, Shi and Yor [11] made use of Lamperti's representation combined with the study of Jacobi processes. Unfortunately, both methods fail for more general potentials. Our approach consists in reducing the study of $H(r)$ to that of an additive functional of a Markov process.

More precisely, the remainder of this paper is organized as follows: in Section 2 we show that $H(r)$ has the same rates of convergence as $\int_0^r Z_s \, ds$, where $Z$ is a generalized Ornstein–Uhlenbeck process. In Section 3 we study the basic properties of $Z$. Section 4 is devoted to the study of the hitting times of $Z$ and we will prove that this process is recurrent. In Section 5 we define the local time and excursion measure associated with the excursions of $Z$ away from level 1. The main result of that section is an estimate of the distribution tail of the area of a generic excursion. Section 6 is devoted to the calculus of the second moment of the area of an excursion when $\kappa > 2$. Once all these results are obtained, the rest of the proof is very classical and is given in the last section.

In the rest of the paper, given a stochastic process $\zeta$, we will indifferently write $\zeta_t$ or $\zeta(t)$.

**2. The process $Z$.** We first construct $X$ from a Brownian motion through a random change of time and a random change of scale. Let $B$ denote a standard Brownian motion independent of $\mathbb{V}$ and for $x \in \mathbb{R}$, set $\sigma_B(x) \stackrel{\text{def}}{=} \inf\{t \geq 0, B_t = x\}$. Recall that the scale function $A$ was defined in (1.3). The process $A$ is continuous and strictly increasing. Let $A^{-1} : (-\infty, A(+\infty)) \mapsto \mathbb{R}$ denote the inverse of $A$. We also define

$$(2.1) \quad T(t) \stackrel{\text{def}}{=} \int_0^t \exp(-2\mathbb{V}(A^{-1}(B_s)))\, ds \qquad \text{for } 0 \leq t < \sigma_B(A(+\infty)).$$

The process $T$ is strictly increasing on $[0, \sigma_B(A(+\infty)))$. Let $T^{-1}$ denote the inverse of $T$ and set

$$(2.2) \qquad\qquad X_t = A^{-1}(B(T^{-1}(t))) \qquad \text{for all } t \geq 0.$$

According to Brox [4], the process $(X_t, t \geq 0)$ is a diffusion in the random potential $\mathbb{V}$. Recall that $H(r)$ defined by (1.6) stands for the hitting time of level $r$ for $X$. Using the representation (2.2), we obtain

$$(2.3) \qquad\qquad H(r) = T(\sigma_B(A(r))).$$

Now, let $L_B(x, t)$ denote the (bi-continuous) local time of $B$ at level $x \in \mathbb{R}$ and time $t \geq 0$. Substituting (2.1) in (2.3), we get

$$H(r) = \int_0^{\sigma_B(A(r))} \exp(-2\mathbb{V}(A^{-1}(B_s)))\, ds$$

$$= \int_{-\infty}^{A(r)} \exp(-2\mathbb{V}(A^{-1}(y))) L_B(y, \sigma_B(A(r)))\, dy.$$

Making use of the change of variable $A(x) = y$,

$$H(r) = \int_{-\infty}^r \exp(-\mathbb{V}_x) L_B(A(x), \sigma_B(A(r)))\, dx = J_1(r) + J_2(r),$$

where

$$J_1(r) \stackrel{\text{def}}{=} \int_{-\infty}^0 \exp(-\mathbb{V}_x) L_B(A(x), \sigma_B(A(r)))\, dx,$$

$$J_2(r) \stackrel{\text{def}}{=} \int_0^r \exp(-\mathbb{V}_x) L_B(A(x), \sigma_B(A(r)))\, dx.$$

We first deal with $J_1$. Since $x \mapsto L_B(x, t)$ has compact support for all $t$ and since $\lim_{x \to -\infty} A(x) = -\infty$, we see that

$$J_1(\infty) \stackrel{\text{def}}{=} \int_{-\infty}^0 \exp(-\mathbb{V}_x) L_B(A(x), \sigma_B(A(+\infty)))\, dx < \infty, \qquad \mathbf{P}\text{-a.s.}$$

Moreover, $J_1(r) \leq J_1(\infty)$ for all $r \geq 0$. Thus, we only need to prove Theorem 1.1 for $J_2(r)$ in place of $H(r)$.

According to the first Ray–Knight theorem, for all $a > 0$, the process $(L_B(a-t, \sigma(a)), 0 \leq t \leq a)$ has the law of a two-dimensional squared Bessel process starting from 0 and is independent of $\mathbb{V}$. Let $(U(x), x \geq 0)$ under $\mathbf{P}$ be a two-dimensional squared Bessel process starting from 0, independent of $\mathbb{V}$. Then, for each fixed $r > 0$,

$$J_2(r) \overset{\text{law}}{=} \int_0^r e^{-\mathbb{V}_x} U(A(r) - A(x)) \, dx$$

$$\overset{\text{law}}{=} \int_0^r e^{-\mathbb{V}_{r-y}} U(A(r) - A(r-y)) \, dy$$

$$\overset{\text{law}}{=} \int_0^r e^{-\mathbb{V}_{r-y}} U\left(\int_0^y e^{\mathbb{V}_{r-s}} \, ds\right) dy$$

$$\overset{\text{law}}{=} \int_0^r e^{-\mathbb{V}_{(r-y)-}} U\left(\int_0^y e^{\mathbb{V}_{(r-s)-}} \, ds\right) dy$$

(where $\mathbb{V}_{x-}$ denotes the left limit of $\mathbb{V}$ at point $x$). For any fixed $r > 0$, we define $\widehat{\mathbb{V}}_t^r \overset{\text{def}}{=} \mathbb{V}_{(r-t)-} - \mathbb{V}_r$ for all $0 \leq t \leq r$. Therefore, the scaling property of $U$ yields

$$J_2(r) \overset{\text{law}}{=} \int_0^r e^{-\widehat{\mathbb{V}}_y^r - \mathbb{V}_r} U\left(e^{\mathbb{V}_r} \int_0^y e^{\widehat{\mathbb{V}}_s^r} \, ds\right) dy$$

$$\overset{\text{law}}{=} \int_0^r e^{-\widehat{\mathbb{V}}_y^r} U\left(\int_0^y e^{\widehat{\mathbb{V}}_s^r} \, ds\right) dy.$$

Time reversal of the Lévy process $\mathbb{V}$ (see Lemma 2, page 45 of [1]) states that, for each $r > 0$, the two processes $(\widehat{\mathbb{V}}_t^r, 0 \leq t \leq r)$ and $(-\mathbb{V}_t, 0 \leq t \leq r)$ have the same law. Thus, for each fixed $r$, under $\mathbf{P}$,

$$(2.4) \qquad J_2(r) \overset{\text{law}}{=} \int_0^r e^{\mathbb{V}_y} U\left(\int_0^y e^{-\mathbb{V}_s} \, ds\right) dy = \int_0^r Z_s \, ds,$$

with the notation

$$(2.5) \qquad\qquad Z_t \overset{\text{def}}{=} e^{\mathbb{V}_t} U(a(t)),$$

and where

$$(2.6) \qquad\qquad a(x) \overset{\text{def}}{=} \int_0^x e^{-\mathbb{V}_s} \, ds.$$

According to (2.4), we only need to prove Theorem 1.1 for the additive functional $\int_0^r Z_s \, ds$ instead of dealing directly with $H(r)$.

The rest of the proof now relies on the study of the process $Z$. As we will see in the next sections, $Z$ is a "nice" recurrent Markov process for which we may define a local time $L$ at any positive level, say, 1. We may therefore also consider the associated excursion measure **n** of its excursions away from 1. Given a generic excursion $(\epsilon_t)$ with lifetime $\zeta$, we define the functional

$$\widetilde{I}(\epsilon) \stackrel{\text{def}}{=} \int_0^\zeta \epsilon_s \, ds.$$

The key step consists in proving that $\widetilde{I}(\epsilon)$, under the excursion measure **n**, has a regularly varying tail of the form

$$\mathbf{n}\{\widetilde{I}(\epsilon) > x\} \underset{x \to \infty}{\sim} \frac{C}{x^\kappa}.$$

Then, as we may write

$$\int_0^t Z_s \, ds \approx \sum_{\substack{\text{excursion } \epsilon \\ \text{starting before } t}} \widetilde{I}(\epsilon),$$

the asymptotics of $\int_0^t Z_s \, ds$ will follow from classical results on the characterization of the domains of attraction to a stable law.

**3. Basic properties of $Z$.**  Recall that $U$ under $\mathbf{P}$ is a two-dimensional squared Bessel process starting from 0 and is independent of $\mathbb{V}$. We now consider a family of probabilities $(\mathbf{P}_x, x \geq 0)$ such that $U$ under $\mathbf{P}_x$ is a two-dimensional squared Bessel process starting from $x$ and is independent of $\mathbb{V}$. In particular, $\mathbf{P} = \mathbf{P}_0$. We will use the notation $\mathbf{E}_x$ for the expectation under $\mathbf{P}_x$ (and $\mathbf{E} = \mathbf{E}_0$ for the expectation under $\mathbf{P} = \mathbf{P}_0$). Of course, the law of $\mathbb{V}$ is the same under all $\mathbf{P}_x$ and when dealing with probabilities that do not depend on the starting point $x$ of $U$, we will use the notation $\mathbf{P}$.

Let us first notice that the process $Z$ defined by (2.5) is nonnegative and does not possess positive jumps because $\mathbb{V}$ has no positive jumps. Moreover, under $\mathbf{P}_x$, the process $Z$ starts from $x$. We define the filtration

$$\mathcal{F}_t \stackrel{\text{def}}{=} \sigma(\mathbb{V}_s, U(a(s)), s \leq t).$$

Our first lemma states that $Z$ is an $\mathcal{F}$-Markov process.

LEMMA 3.1.  *$((Z_t)_{t \geq 0}, (\mathbf{P}_x)_{x \geq 0})$ is an $\mathcal{F}$-Markov process whose semigroup fulfills the Feller property. Moreover, for each $x > 0$, the process $(Z_t, t \geq 0)$ under $\mathbf{P}_x$ (i.e., starting from $x$) has the same law as the process $(\widetilde{Z}_t^x, t \geq 0)$ under $\mathbf{P}_1$, where*

$$(3.1) \qquad \widetilde{Z}_t^x \stackrel{\text{def}}{=} x e^{\mathbb{V}_t} U\left(\frac{a(t)}{x}\right).$$

Proof. The process $U$ is a squared Bessel process. Therefore, our process $Z$ is a generalized Ornstein–Uhlenbeck process in the sense of [6] and Proposition 5.5 of [6] states that $Z$ is indeed a Markov process in the filtration $\mathcal{F}$. Let $(P_t)_{t \geq 0}$ and $(Q_t)_{t \geq 0}$ stand for the respective semi-groups of $U$ and $Z$. The independence of $U$ and $\mathbb{V}$ yields the relation

$$(3.2) \qquad Q_t f(x) = \mathbf{E}_x[f(Z_t)] = \mathbf{E}[P_{a(t)}(f(e^{\mathbb{V}_t} \cdot))(x)].$$

Since $U$ is a squared Bessel process, its semi-group fulfills the Feller property. Moreover, $a(\cdot)$ is continuous with $a(0) = 0$ and $\lim_{t \to 0+} e^{\mathbb{V}_t} = 1$ $\mathbf{P}$-a.s. These facts combined with (3.2) easily show that $(Q_t)$ is also a Fellerian semi-group. Finally, (3.1) is an immediate consequence of the scaling property of $U$. □

For $x \geq 0$, we say that $x$ is instantaneous for $Z$ if the process $Z$ starting from $x$ leaves $x$ instantaneously with probability 1. Moreover, we say $x$ is regular (for itself) for $Z$ if $Z$ starting from $x$ returns to $x$ at arbitrarily small times with probability 1.

Lemma 3.2. *Any $x > 0$ is regular and instantaneous for $Z$.*

Proof. We only prove the result for $x = 1$; the general case may be treated the same way. Since $U$ under $\mathbf{P}_1$ is a squared Bessel process of dimension 2 starting from 1, it has the same law as $(B^2(t) + \widetilde{B}^2(t) + 2B(t) + 1, t \geq 0)$, where $B$ and $\widetilde{B}$ are two independent standard Brownian motions. It is therefore easy to check using classical results on Brownian motion that:

(a) For any strictly decreasing sequence $(t_i)_{i \geq 0}$ of (nonrandom) real numbers with $\lim_{i \to \infty} t_i = 0$, we have

$$\mathbf{P}_1\{U(t_i) > 1 \text{ i.o.}\} = \mathbf{P}_1\{U(t_i) < 1 \text{ i.o.}\} = 1.$$

(b) $\liminf_{t \to 0+} \frac{U(t) - 1}{t} = -\infty$, $\mathbf{P}_1$-a.s.

Let us now prove that $Z$ starting from 1 visits $(1, \infty)$ at arbitrarily small times. Recall that $(\mathbb{V}_x, x \geq 0)$ is a Lévy process with no positive jumps which is not the opposite of a subordinator. According to Theorem 1, page 189 of [1], the process $\mathbb{V}$ visits $(0, \infty)$ at arbitrarily small times with probability 1. Thus, for almost any fixed path of $\mathbb{V}$, we can find a strictly positive decreasing sequence $(u_i)_{i \geq 0}$ with limit 0 such that $\mathbb{V}_{u_i} > 0$ for all $i$. But, conditionally on $\mathbb{V}$, under $\mathbf{P}_1$, $U$ is still a squared Bessel process of dimension 2 starting from 1 and

$$Z_{u_i} = e^{\mathbb{V}_{u_i}} U(a(u_i)) > U(a(u_i)).$$

Since $a(\cdot)$ is continuous with $\lim_{t\to 0} a(t) = 0$, the sequence $(a(u_i))_{i\geq 0}$ is positive, strictly decreasing with limit 0. Using (a), we conclude that $Z$ starting from 1 visits $(1,\infty)$ at arbitrarily small times almost surely.

When 0 is regular for $(-\infty, 0)$ for the Lévy process $\mathbb{V}$, a similar argument shows that $Z$ starting from 1 visits $(0,1)$ at arbitrarily small times almost surely. Let us therefore assume that 0 is irregular for $(-\infty, 0)$ for $\mathbb{V}$. According to Corollary 5, page 192 of [1], this implies that $\mathbb{V}$ has bounded variations, thus, there exists $d \geq 0$ such that $\lim_{x\to 0+} \mathbb{V}_x/x = d$ a.s. (cf. Proposition 11, page 166 of [1]). Let $a^{-1}(\cdot)$ denote the inverse of $a(\cdot)$. Since $a(t) \sim t$ as $t \to 0+$, we have $e^{\mathbb{V}_{a^{-1}(t)}} \leq 1 + 2\,dt$ for all $t$ small enough, almost surely. In consequence,

$$Z(a^{-1}(t)) = e^{\mathbb{V}_{a^{-1}(t)}}U(t) \leq (1+2\,dt)U(t) \qquad \text{for } t \text{ small enough, } \mathbf{P}_1\text{-a.s.}$$

Using (b), we conclude that the process $(Z(a^{-1}(t)), t \geq 0)$ visits $(0,1)$ at arbitrarily small times $\mathbf{P}_1$-a.s. Since $a^{-1}(\cdot)$ is continuous, increasing and $a(0) = 0$, this result also holds for Z.

We proved that $Z$ starting from 1 visits $(0,1)$ and $(1,\infty)$ at arbitrarily small times almost surely. Since $Z$ has no positive jumps, $Z$ starting from 1 returns to 1 at arbitrarily small times almost surely.  $\square$

LEMMA 3.3.   *For all $x, y \geq 0$ and all $t > 0$, we have $\mathbf{P}_x\{Z_t = y\} = 0$. In consequence,*

$$\int_0^\infty \mathbf{1}_{\{Z_t = y\}}\,dt = 0, \qquad \mathbf{P}_x\text{-a.s. for all } x, y \geq 0.$$

PROOF.   A squared Bessel process has a continuous density, in particular, $\mathbf{P}_x\{U(a) = b\} = 0$ for all $b, x \geq 0$ and all $a > 0$. Since $\mathbb{V}$ and $U$ are independent and $a(t) > 0$ for all $t > 0$, we get

$$\mathbf{P}_x\{Z_t = y\} = \mathbf{E}[\mathbf{P}_x\{U(a(t)) = ye^{-\mathbb{V}_t}|\mathbb{V}\}] = 0. \qquad \square$$

The following easy lemma will be found very useful in the remainder of this paper.

LEMMA 3.4.   *For all $0 \leq x \leq y$, the process $Z$ under $\mathbf{P}_x$ (i.e., starting from $x$) is stochastically dominated by $Z$ under $\mathbf{P}_y$ (i.e., starting from $y$).*

PROOF.   The process $U$ is a two-dimensional squared Bessel process and a theorem of comparison for diffusion process (cf. Theorem IX.3.7 of [16]) states that $U$ under $\mathbf{P}_x$ is stochastically dominated by $U$ under $\mathbf{P}_y$ whenever $x \leq y$. Thus, the lemma is a direct consequence of the independence of $U$ and $\mathbb{V}$.  $\square$

We conclude this section by proving the convergence of $Z$ at infinity.

PROPOSITION 3.1. *Letting $x > 0$, under $\mathbf{P}_x$, $Z_t$ converges as $t$ goes to infinity toward a nondegenerate random variable $Z_\infty$ whose law does not depend on the starting point $x$. The distribution of $Z_\infty$ is the same as that of the random variable*

$$(3.3) \qquad U(1)\int_0^\infty e^{\mathbb{V}_s}\, ds \qquad under \ \mathbf{P}_0.$$

*In particular, the law of $Z_\infty$ has a strictly positive continuous density on $(0, \infty)$ and*

$$(3.4) \qquad \mathbf{P}\{Z_\infty > x\} \underset{x \to \infty}{\sim} \frac{2^\kappa \Gamma(\kappa+1)\mathbf{K}}{\Phi'(\kappa)x^\kappa},$$

*where $\kappa$ is the constant of (1.2) and where $\mathbf{K} = \mathbf{E}[A(+\infty)^{\kappa-1}] \in (0, \infty)$ is the constant defined in the statement of Theorem 1.1.*

PROOF. According to Proposition 5.7 of [6], under the assumption that $\mathbf{E}[\mathbb{V}_1] < 0$, the generalized Ornstein–Uhlenbeck process $Z$ converges in law toward a random variable $Z_\infty$ whose distribution is given by (3.3). In our case, we may also have $\mathbf{E}[\mathbb{V}_1] = -\infty$. However, in the proof of Proposition 5.7 of [6], the assumption that $E[\mathbb{V}_1] < 0$ is required only to ensure that

$$\lim_{t \to \infty} \mathbb{V}_t = -\infty \quad \text{and} \quad \int_0^\infty e^{\mathbb{V}_t}\, dt = A(+\infty) < \infty \qquad \text{a.s.}$$

Since we have already established these two results, Proposition 5.7 of [6] is also true in our case. The process $U$ under $\mathbf{P}_0$ is a squared Bessel process of dimension 2 starting from 0, therefore, $U(1)$ under $\mathbf{P}_0$ has an exponential distribution with mean 2. Keeping in mind that $\mathbb{V}$ and $U$ are independent, we find

$$
\begin{aligned}
\mathbf{P}\{Z_\infty > x\} &= \mathbf{P}_0\{U(1)A(+\infty) > x\} \\
&= \mathbf{E}[\mathbf{P}_0\{U(1)A(+\infty) > x | A(+\infty)\}] \\
(3.5) \qquad &= \mathbf{E}\Big[\exp\Big(-\frac{x}{2A(+\infty)}\Big)\Big].
\end{aligned}
$$

It is now clear that $Z_\infty$ has a continuous density, everywhere positive on $(0, \infty)$. Moreover, in view of the Abelian/Tauberian theorem (see, e.g., Chapter VIII of [9]), we deduce from (3.5) that the estimate (3.4) on the tail distribution of $Z_\infty$ is equivalent to

$$(3.6) \qquad \mathbf{P}\{A(+\infty) > x\} \underset{x \to \infty}{\sim} \frac{\mathbf{E}[A(+\infty)^{\kappa-1}]}{\Phi'(\kappa)x^\kappa}.$$

This result is proved in Lemma 4 of [17] in the case $0 < \kappa < 1$. Another proof, valid for any $\kappa > 0$, is given in Theorem 3.1 of [15] under the restrictive assumption that $\mathbb{V}_1$ admits a finite first moment. However, one may check that, in the proof of Theorem 3.1 of [15], the assumption $\mathbf{E}[|\mathbb{V}_1|] < \infty$ is only needed for $0 < \kappa < 1$. Thus, in our setting, (3.6) holds for any $\kappa > 0$. We point out that Lemma 4 of [17] and Theorem 3.1 of [15] are both based on a theorem of Goldie [10] which is, in turn, a refined version in the one-dimensional case of a famous result of Kesten [13] on the affine equation for random matrices. □

**4. Hitting times of $Z$.** Given a stochastic process $Y$ and a set $A$, we define the hitting times

(4.1)    $\tau_A(Y) = \inf\{t \geq 0, Y_t \in A\}$      (with the convention $\inf \varnothing = \infty$).

For simplicity, we will use the notation $\tau_x(Y)$ instead of $\tau_{\{x\}}(Y)$. When referring to the process $Z$, we will also simply write $\tau_A$ instead of $\tau_A(Z)$. We now show that the hitting times of $Z$ are finite almost-surely and we give estimates on their distribution tail. In particular, this will show that $Z$ is recurrent. The rest of this section is devoted to proving the following four propositions. These estimates are quite technical and, on a first reading, the details of the proof may be skipped over after glancing at the statements of the propositions.

PROPOSITION 4.1.    *For any $0 \leq x < y$, there exist $c_{1,y}, c_{2,y} > 0$ (depending on $y$) such that*

$$\mathbf{P}_x\{\tau_{[y,\infty)} > t\} \leq c_{1,y} e^{-c_{2,y}t}      \text{for all } t \geq 0.$$

PROPOSITION 4.2.    *There exist $y_0, c_3, c_4 > 0$ such that, for all $y_0 \leq y < x$,*

$$\mathbf{P}_x\{\tau_{[0,y]} > t\} \leq c_3(\log(x/y) + 1)e^{-c_4/(\log(x/y)+1)t}      \text{for all } t \geq 0.$$

PROPOSITION 4.3.    *For all $x \geq 0$ and all $y > 0$, there exist $c_{5,x,y}, c_{6,x,y} > 0$ (depending on $x$ and $y$) such that*

$$\mathbf{P}_x\{\tau_y > t\} \leq c_{5,x,y} e^{-c_{6,x,y}t}      \text{for all } t \geq 0.$$

*In particular, $Z$ starting from $x \geq 0$ hits any positive level eventually.*

PROPOSITION 4.4.    *We have*

$$\lim_{\lambda \to \infty} \sup_{y \geq 1} \mathbf{P}_y\{\tau_{\lambda y} < \tau_1\} = 0.$$

PROOF OF PROPOSITION 4.1.    Let $0 \leq x \leq y$. According to Lemma 3.4, $\tau_{[y,\infty)}$ under $\mathbf{P}_0$ is stochastically dominated by $\tau_{[y,\infty)}$ under $\mathbf{P}_x$, thus, we only need to prove the proposition for $x = 0$. Let $\lfloor t \rfloor$ stand for the integer part of $t$. We have

$$\mathbf{P}_0\{\tau_{[y,\infty)} > t\} \leq \mathbf{P}_0\{Z_1 < y, Z_2 < y, \ldots, Z_{\lfloor t \rfloor} < y\} \leq \mathbf{P}_0\{Z_1 < y\}^{\lfloor t \rfloor},$$

where we repeatedly used the Markov property of $Z$ combined with the stochastic monotonicity of $Z$ (Lemma 3.4) for the last inequality. Since $Z_1 = e^{\mathbb{V}_1}U(a(1))$, it is clear that $\mathbf{P}_0\{Z_1 < y\} < 1$ for all $y > 0$. Thus, setting $c_{2,y} = -\log(\mathbf{P}_0\{Z_1 < y\}) > 0$ and $c_{1,y} = e^{c_{2,y}}$, we find

$$\mathbf{P}_0\{\tau_{[y,\infty)} > y\} \leq e^{-c_{2,y}\lfloor t \rfloor} \leq c_{1,y}e^{-c_{2,y}t}. \qquad \square$$

The proof of Proposition 4.2 relies on the following:

LEMMA 4.1.    *There exist $c_7, c_8, x_0 > 0$ such that, for all $x \geq x_0$,*

$$\mathbf{P}_x\{\tau_{[0,x/2]} > t\} \leq c_7 e^{-c_8 t} \qquad \text{for all } t > 0.$$

PROOF.    Pick $\eta > 0$ and let $(\mathbb{V}_t^{(\eta)}, t \geq 0)$ stand for the Lévy process $\mathbb{V}_t^{(\eta)} = \mathbb{V}_t + \eta t$. Recall that $\Phi$ denotes the Laplace exponent of $\mathbb{V}$. Thus, the Laplace exponent $\Phi^{(\eta)}$ of $\mathbb{V}^{(\eta)}$ is given by $\Phi^{(\eta)}(x) = \Phi(x) + \eta x$. Since $\Phi(\kappa/2) < 0$, we can choose $\eta$ small enough such that $\Phi^{(\eta)}(\kappa/2) < 0$. Then $\mathbb{V}_t^{(\eta)}$ diverges to $-\infty$ as $t$ goes to infinity and we can define the sequence

$$\begin{cases} \gamma_0 \overset{\text{def}}{=} 0, \\ \gamma_{n+1} \overset{\text{def}}{=} \inf\{t > \gamma_n, \mathbb{V}_t^{(\eta)} - \mathbb{V}_{\gamma_n}^{(\eta)} < -\log(8)\}. \end{cases}$$

The sequence $(\gamma_{n+1} - \gamma_n)_{n \geq 0}$ is i.i.d. and distributed as $\gamma_1$. We have

$$\mathbf{P}\{\gamma_1 > t\} \leq \mathbf{P}\{\mathbb{V}_t^{(\eta)} \geq -\log(8)\} \leq \mathbf{P}\left\{\exp\left(\frac{\kappa}{2}\mathbb{V}_t^{(\eta)}\right) \geq \frac{1}{8^{\kappa/2}}\right\}$$

$$\leq 8^{\kappa/2}\mathbf{E}\left[\exp\left(\frac{\kappa}{2}\mathbb{V}_t^{(\eta)}\right)\right] = 8^{\kappa/2}e^{t\Phi^{(\eta)}(\kappa/2)}.$$

Since $\Phi^{(\eta)}(\kappa/2) < 0$, we deduce from Cramér's large deviation theorem that there exist $c_9, c_{10}, c_{11} > 0$ such that

$$(4.2) \qquad \mathbf{P}\{\gamma_n > c_9 n\} \leq c_{10}e^{-c_{11}n} \qquad \text{for all } n \in \mathbb{N}.$$

Notice from the definition of $\gamma_1$ that

$$(4.3) \quad e^{\mathbb{V}_{\gamma_1}}a(\gamma_1) = \int_0^{\gamma_1} e^{\mathbb{V}_{\gamma_1}^{(\eta)} - \mathbb{V}_s^{(\eta)} - \eta(\gamma_1 - s)}\, ds \leq \int_0^{\gamma_1} e^{-\eta(\gamma_1 - s)}\, ds \leq \frac{1}{\eta},$$

and also

$$(4.4) \qquad\qquad\qquad e^{\mathbb{V}_{\gamma_1}} \leq \tfrac{1}{8}.$$

The process $U$ under $\mathbf{P}_x$ is a squared Bessel process of dimension 2 starting from $x$. Therefore, $U$ under $\mathbf{P}_x$ is stochastically dominated by $2(x+U)$ under $\mathbf{P}_0$. Using the independence of $\mathbb{V}$ and $U$, we deduce that $Z_{\gamma_1}$ under $\mathbf{P}_x$ is stochastically dominated by $2e^{\mathbb{V}_{\gamma_1}}(x+U(a(\gamma_1)))$ under $\mathbf{P}_0$. Moreover, the scaling property of $U$ combined with (4.3) and (4.4) yields, under $\mathbf{P}_0$,

$$2e^{\mathbb{V}_{\gamma_1}}\left(x+U(a(\gamma_1))\right) \overset{\text{law}}{=} 2xe^{\mathbb{V}_{\gamma_1}}+2e^{\mathbb{V}_{\gamma_1}}a(\gamma_1)U(1) \le \frac{x}{4}+\frac{2}{\eta}U(1).$$

Thus, $Z_{\gamma_1}$ under $\mathbf{P}_x$ is stochastically dominated by the random variable $\frac{x}{4}+\frac{2}{\eta}U(1)$ under $\mathbf{P}_0$. Now, let $(\chi_n, n \ge 1)$ denote a sequence of i.i.d. random variables with the same distribution as $\frac{2}{\eta}U(1)$ under $\mathbf{P}_0$. Define also the sequence $(R_n^x, n \ge 0)$ by

$$\begin{cases} R_0^x \overset{\text{def}}{=} x, \\ R_{n+1}^x \overset{\text{def}}{=} \frac{1}{4}R_n^x + \chi_{n+1}. \end{cases}$$

The process $(Z_{\gamma_n}, n \ge 0)$ under $\mathbf{P}_x$ is a Markov chain starting from $x$. We already proved that $Z_{\gamma_1}$ is stochastically dominated by $R_1^x$. By induction and with the help of Lemma 3.4, we conclude with similar arguments that the sequence $(Z_{\gamma_n}, n \ge 0)$ under $\mathbf{P}_x$ is stochastically dominated by $(R_n^x, n \ge 0)$. In particular, choosing $n = \lfloor t/c_9 \rfloor$ and using (4.2), we find

$$\mathbf{P}_x\{\tau_{[0,x/2]} > t\} \le \mathbf{P}\{\gamma_n > c_9 n\} + \mathbf{P}_x\left\{ Z_{\gamma_1} > \frac{x}{2}, \ldots, Z_{\gamma_n} > \frac{x}{2} \right\}$$

$$\le c_{10}e^{-c_{11}t} + \mathbf{P}\left\{ R_1^x > \frac{x}{2}, \ldots, R_n^x > \frac{x}{2} \right\}.$$

Thus, it only remains to prove that there exist $c_{12}, x_0 > 0$ such that

$$\mathbf{P}\left\{ R_1^x > \frac{x}{2}, \ldots, R_n^x > \frac{x}{2} \right\} \le e^{-c_{12}n} \qquad \text{for all } n \in \mathbb{N} \text{ and all } x \ge x_0.$$

Expanding the definition of $R^x$, we get

$$R_n^x = \frac{x}{4^n} + \frac{1}{4^{n-1}}\chi_1 + \cdots + \frac{1}{4}\chi_{n-1} + \chi_n.$$

Let us set $c = 8/\eta$. We have

$$(4.5) \quad R_n^x - \frac{4}{3}c \le \frac{x}{4^n} + \frac{1}{4^{n-1}}(\chi_1 - c) + \cdots + \frac{1}{4}(\chi_{n-1}-c) + (\chi_n - c).$$

Let also $S$ denote the random walk given by $S_0 \overset{\text{def}}{=} 0$ and $S_{n+1} \overset{\text{def}}{=} S_n + (\chi_{n+1}-c)$. We can rewrite (4.5) in the form

$$(4.6) \quad R_n^x - \frac{4}{3}c \le \frac{x}{4^n} + S_n - \frac{3}{4}\sum_{k=1}^{n-1}\frac{1}{4^{n-1-k}}S_k.$$

Let $\mu \overset{\text{def}}{=} \inf\{n \geq 1, S_n < 0\}$ stand for the first strict descending ladder index of the random walk $S$. We have

$$\mathbf{P}\{\mu > n\} \leq \mathbf{P}\{S_n \geq 0\} \leq \mathbf{E}[e^{(\eta/8)S_n}] = \mathbf{E}[e^{(\eta/8)S_1}]^n = \left(\frac{2}{e}\right)^n,$$

where we used the fact that $S_1$ has the same distribution as the random variable $\frac{2}{\eta}U(1) - \frac{8}{\eta}$ under $\mathbf{P}_0$ [and $U(1)$ under $\mathbf{P}_0$ has an exponential distribution with mean 2]. Therefore, $\mu$ is almost surely finite. Setting $c_{12} = 1 - \log 2 > 0$, we get

$$\mathbf{P}\{\mu > n\} \leq e^{-c_{12}n} \qquad \text{for all } n \in \mathbb{N}.$$

Finally, from the definition of $\mu$, we have $S_\mu < 0$ and $S_n \geq 0$ for all $0 \leq n < \mu$, hence,

$$S_\mu - \frac{3}{4}\sum_{k=1}^{\mu-1}\frac{1}{4^{n-1-k}}S_k \leq 0.$$

Combining this inequality with (4.6) and the fact that $\mu \geq 1$, we obtain, whenever $x \geq x_0 \overset{\text{def}}{=} \frac{16}{3}c$

$$R_\mu^x \leq \frac{4}{3}c + \frac{x}{4^\mu} \leq \frac{4}{3}c + \frac{x}{4} \leq \frac{x}{2}.$$

Thus, for all $x \geq x_0$,

$$\mathbf{P}\left\{R_1^x > \frac{x}{2}, \ldots, R_n^x > \frac{x}{2}\right\} \leq \mathbf{P}\{\mu > n\} \leq e^{-c_{12}n}.$$

This completes the proof of the lemma.   $\square$

PROOF OF PROPOSITION 4.2.   Set $y_0 \overset{\text{def}}{=} x_0/2$, where $x_0$ is the constant of the previous lemma. This lemma ensures that for all $x, y$ such that $y_0 < y < x$ and $\frac{x}{y} \leq 2$, we have

$$(4.7) \qquad \mathbf{P}_x\{\tau_{[0,y]} > t\} \leq c_7 e^{-c_8 t} \qquad \text{for all } t > 0.$$

Let us now fix $x, y$ such that $y_0 \leq y < x$. Define the sequence $(z_n)$ by $z_0 \overset{\text{def}}{=} x$ and $z_{n+1} \overset{\text{def}}{=} z_n/2$. We also set $m \overset{\text{def}}{=} 1 + \lfloor \log(x/y)/\log(2)\rfloor$, then

$$x = z_0 \geq z_1 \geq \cdots \geq z_{m-1} \geq y \geq z_m.$$

Thus,

$$\mathbf{P}_x\{\tau_{[0,y]} > t\} \leq \mathbf{P}_x\left\{\tau_{[0,z_1]} > \frac{t}{m}\right\} + \sum_{i=1}^{m-2}\mathbf{P}_x\left\{\tau_{[0,z_{i+1}]} - \tau_{[0,z_i]} > \frac{t}{m}\right\}$$

$$+ \mathbf{P}_x\left\{\tau_{[0,y]} - \tau_{[0,z_{m-1}]} > \frac{t}{m}\right\}.$$

Making use of the Markov property of $Z$ for the stopping times $\tau_{[0,x_i]}$ combined with Lemma 3.4, we get

$$\mathbf{P}_x\{\tau_{[0,y]} > t\} \leq \mathbf{P}_{z_0}\left\{\tau_{[0,z_1]} > \frac{t}{m}\right\} + \sum_{i=1}^{m-2} \mathbf{P}_{z_i}\left\{\tau_{[0,z_{i+1}]} > \frac{t}{m}\right\}$$
$$+ \mathbf{P}_{z_{m-1}}\left\{\tau_{[0,y]} > \frac{t}{m}\right\}.$$

According to (4.7), each term on the r.h.s. of this last inequality is smaller than $c_7 e^{-c_8 t/m}$, hence, choosing $c_3, c_4$ large enough,

$$\mathbf{P}_x\{\tau_{[0,y]} > t\} \leq mc_7 e^{-c_8 t/m} \leq c_3(\log(x/y) + 1)e^{-c_4(\log(x/y)+1)t}. \qquad \square$$

PROOF OF PROPOSITION 4.3. We have already proved Proposition 4.1 and Proposition 4.2. Recall also that, $Z$ has no positive jumps. In view of Lemma 3.4, it simply remains to prove that for any $0 < y < y_0$ ($y_0$ is the constant of Proposition 4.2), we have

(4.8) $$\mathbf{P}_{y_0}\{\tau_{[0,y]} > t\} \leq c_{13,y_0,y} e^{-c_{14,y_0,y}t} \qquad \text{for all } t > 0.$$

Let us fix $y < y_0$. We also pick $z > y_0$. Define the sequence $(\nu_n^z)$:

$$\begin{cases} \nu_0^z \overset{\text{def}}{=} 0, \\ \nu_{n+1}^z \overset{\text{def}}{=} \inf\left\{t > \nu_n^z, Z_t = y_0 \text{ and } \sup_{\nu_n^z \leq s \leq t} Z_s \geq z\right\}. \end{cases}$$

Making use of Propositions 4.1 and 4.2, we check that $\nu_n^z$ is finite for all $n$, $\mathbf{P}_{y_0}$-a.s. More precisely, these propositions yield

$$\mathbf{P}_{y_0}\{\nu_1^z > t\} \leq c_{15,y_0,z} e^{-c_{16,y_0,z}t} \qquad \text{for all } t > 0.$$

Since the sequence $(\nu_{n+1}^z - \nu_n^z)_{n\geq 0}$ is i.i.d., Cramér's large deviation theorem ensures that there exist $c_{17,y_0,z}, c_{18,y_0,z}, c_{19,y_0,z} > 0$ such that

(4.9) $$\mathbf{P}_{y_0}\{\nu_n^z > c_{17,y_0,z}n\} \leq c_{18,y_0,z} e^{-c_{19,y_0,z}n} \qquad \text{for all } n \in \mathbb{N}.$$

Let us note that $\lim_{z\to\infty} \nu_1^z = \infty$, $\mathbf{P}_{y_0}$-a.s., thus,

(4.10) $$\mathbf{P}_{y_0}\{\tau_{[0,y]} < \nu_1^z\} \underset{z\to\infty}{\longrightarrow} \mathbf{P}_{y_0}\{\tau_{[0,y]} < \infty\}.$$

According to Proposition 3.1, we have $\mathbf{P}_{y_0}\{Z_\infty \in (0,y]\} > 0$. In particular, the limit in (4.10) is strictly positive. Thus, we may choose $z$ large enough such that $\mathbf{P}_{y_0}\{\tau_{[0,y]} > \nu_1^z\} = d < 1$. Repeated use of the Markov property of $Z$ for the stopping times $\nu_i^z$ yields

(4.11) $$\mathbf{P}_{y_0}\{\tau_{[0,y]} > \nu_n^z\} = \mathbf{P}_{y_0}\{\tau_{[0,y]} > \nu_1^z\}^n = d^n.$$

Finally, setting $n = \lfloor t/c_{15,y_0,z} \rfloor$, we get from (4.9) and (4.11)

$$\mathbf{P}_{y_0}\{\tau_{[0,y]} > t\} \leq \mathbf{P}_{y_0}\{\nu_n^z > t\} + \mathbf{P}_{y_0}\{\tau_{[0,y]} > \nu_n^z\}$$
$$\leq c_{18,y_0,z} e^{-c_{19,y_0,z} n} + d^n$$
$$\leq c_{20,y_0,y} e^{-c_{21,y_0,y} t}. \qquad \square$$

We need the following lemma before giving the proof of Proposition 4.4.

Lemma 4.2. *There exist $k_0 > 1$ and $y_1 > 1$ such that*

$$\mathbf{P}_y\{\tau_{k_0 y} < \tau_{y/k_0}\} \leq \tfrac{1}{4} \qquad \text{for all } y \geq y_1.$$

Proof. Let us choose $k > 1$ and $y$ such that

$$(4.12) \qquad 6^4 k^5 < y.$$

We also use the notation $m \overset{\text{def}}{=} \frac{1}{4} \log(\frac{y}{k})$ and $\gamma_m \overset{\text{def}}{=} \inf\{t \geq 0, \mathbb{V}_t < -m\}$. Define $\mathcal{E}_1 \overset{\text{def}}{=} \{\gamma_m \leq e^m\}$. Since $\Phi(\kappa/2) < 0$, we deduce that

$$\mathbf{P}\{\mathcal{E}_1^c\} \leq \mathbf{P}\{\mathbb{V}_{e^m} > -m\} \leq e^{\kappa/2m} \mathbf{E}[e^{\kappa/2 \mathbb{V}_{e^m}}] = e^{(\kappa/2)m + e^m \Phi(\kappa/2)} \underset{y/k \to \infty}{\longrightarrow} 0.$$

We also consider $\mathcal{E}_2 \overset{\text{def}}{=} \{\sup_{s \geq 0} \mathbb{V}_s < \log(k/7)\}$. Since $\mathbb{V}$ diverges to $-\infty$, its overall supremum is finite (it has an exponential distribution with parameter $\kappa$), therefore,

$$\mathbf{P}\{\mathcal{E}_2^c\} \underset{k \to \infty}{\longrightarrow} 0.$$

Define also

$$\mathcal{E}_3 \overset{\text{def}}{=} \left\{ U(t) \leq 2\left(y + \frac{y}{k} + t^2\right) \text{ for all } t \geq 0 \right\}.$$

We noticed in the proof of Lemma 4.1 that $U$ under $\mathbf{P}_y$ is stochastically dominated by $2(y + B^2 + \widetilde{B}^2)$, where $B$ and $\widetilde{B}$ are two independent squared Brownian motions. Therefore, the law of the iterated logarithm for Brownian motion (see, e.g., Chapter II of [16]) entails

$$\mathbf{P}_y\{\mathcal{E}_3^c\} \leq \mathbf{P}_0\left\{ \text{there exists } t \geq 0 \text{ with } U(t) > \frac{y}{k} + t^2 \right\} \underset{y/k \to \infty}{\longrightarrow} 0.$$

We finally set $\mathcal{E}_4 \overset{\text{def}}{=} \mathcal{E}_1 \cap \mathcal{E}_2 \cap \mathcal{E}_3$. Our previous estimates ensure that $\mathbf{P}_y\{\mathcal{E}_4^c\} < 1/4$ whenever $k$ and $y/k$ are both large enough. Moreover, on the set $\mathcal{E}_4$, for all $0 \leq t \leq \gamma_m$,

$$a(t)^2 = \left( \int_0^t e^{-\mathbb{V}_s} \, ds \right)^2 \leq \left( \int_0^{\gamma_m} e^{-\mathbb{V}_s} \, ds \right)^2 \leq (\gamma_m e^m)^2 \leq e^{4m} = \frac{y}{k}.$$

Thus, on the one hand, on $\mathcal{E}_4$, for $k \geq 1$ and for all $0 \leq t \leq \gamma_m$,

$$Z_t = e^{\mathbb{V}_t} U(a(t)) \leq e^{(\sup_{s \geq 0} \mathbb{V}_s)} 2\left(y + \frac{y}{k} + a(t)^2\right) \leq \frac{2k}{7}\left(y + \frac{y}{k} + \frac{y}{k}\right) < ky.$$

On the other hand, on $\mathcal{E}_4$, since $\mathbb{V}_{\gamma_m} \leq -m$,

$$Z_{\gamma_m} \leq e^{-m} 2\left(y + \frac{y}{k} + \frac{y}{k}\right) \leq 6y e^{-m} \leq \frac{y}{k},$$

where we used (4.12) for the last inequality. Therefore,

$$\mathbf{P}_y\{\tau_{[ky,\infty)} < \tau_{[0,y/k]}\} \leq \mathbf{P}_y\{\mathcal{E}_4^c\} < \frac{1}{4} \qquad \text{for all } k, \frac{y}{k} \text{ large enough.}$$

Finally, since $Z$ has no positive jumps, we also have

$$\mathbf{P}_y\{\tau_{[ky,\infty)} < \tau_{[0,y/k]}\} = \mathbb{P}_y\{\tau_{ky} < \tau_{y/k}\}. \qquad \square$$

PROOF OF PROPOSITION 4.4. Let $y_1$ and $k_0$ denote the constants of the previous lemma and let $y \geq y_1$. Define the sequence $(\mu_n)$ of stopping times for $Z$:

$$\begin{cases} \mu_0 \overset{\text{def}}{=} 0, \\ \mu_{n+1} \overset{\text{def}}{=} \inf\left\{t > \mu_n, Z_t = k_0 Z_{\mu_n} \text{ or } Z_t = \frac{1}{k_0} Z_{\mu_n}\right\}. \end{cases}$$

Proposition 4.1 ensures that $\mu_n < \infty$ for all $n$, $\mathbf{P}_y$-a.s. The Markov property of $Z$ also implies that the sequence $(Z_{\mu_n}, n \in \mathbb{N})$ is, under $\mathbf{P}_y$, a Markov chain starting from $y$ and taking values in $\{k_0^n y, n \in \mathbb{Z}\}$. Moreover, according to the previous lemma,

$$\mathbf{P}_y\{Z_{\mu_{n+1}} = k_0 Z_{\mu_n} | Z_{\mu_n} > y_1\} = 1 - \mathbf{P}_y\left\{Z_{\mu_{n+1}} = \frac{1}{k_0} Z_{\mu_n} \Big| Z_{\mu_n} > y_1\right\} < \frac{1}{4}.$$

Thus, if $(S_n, n \geq 0)$ now denotes a random walk such that

$$\begin{cases} \mathbf{P}\{S_0 = 0\} = 1, \\ \mathbf{P}\{S_{n+1} = S_n + 1\} = 1 - \mathbf{P}\{S_{n+1} = S_n - 1\} = \frac{1}{4}, \end{cases}$$

then we deduce from the previous lemma that $(Z_{\mu_n})_{0 \leq n \leq \inf\{n \geq 0, Z_{\mu_n} \leq y_1\}}$ under $\mathbf{P}_y$ is stochastically dominated by $(y k_0^{S_n})_{0 \leq n \leq \inf\{n \geq 0, y k_0^{S_n} \leq y_1\}}$. In particular, for all $y \geq y_1$ and all $p \in \mathbb{N}^*$,

$$\mathbf{P}_y\{(Z_{\mu_n}) \text{ hits } [k_0^p y, \infty) \text{ before it hits } [0, y_1]\} \leq \mathbf{P}\left\{\sup_n S_n \geq p\right\}.$$

Since $Z$ has no positive jumps, we obtain, for all $y \geq y_1$ and all $p \in \mathbb{N}^*$,

$$(4.13) \qquad \mathbf{P}_y\{\tau_{y k_0^p} < \tau_{[0,y_1]}\} \leq \mathbf{P}\left\{\sup_n S_n \geq p\right\}.$$

Note that the last inequality is trivial when $y \leq y_1$. Also, since $S$ is transient toward $-\infty$, its overall supremum is finite and, given $\varepsilon > 0$, we may find $p_0$ such that $\mathbf{P}\{\sup_n S_n \geq p_0\} \leq \varepsilon$. Setting $\lambda_0 \overset{\text{def}}{=} k_0^{p_0}$, we deduce from (4.13) that

$$\sup_{\lambda \geq \lambda_0} \sup_{y \geq 1} \mathbf{P}_y\{\tau_{y\lambda} < \tau_{[0,y_1]}\} \leq \varepsilon.$$

Note that $\tau_{[0,y_1]} \leq \tau_1$ (because $y_1 \geq 1$) and recall that $Z$ has no positive jumps. By use of the Markov property of $Z$ and with the help of Lemma 3.4, we obtain, for all $y \geq 1$ and all $\lambda > \lambda_0$,

$$\mathbf{P}_y\{\tau_{y} < \tau_1\} \leq \mathbf{P}_y\{\tau_{y} < \tau_{[0,y_1]}\} + \mathbf{P}_y\{\tau_{[0,y_1]} \leq \tau_{y}\}\mathbb{P}_{y_1}\{\tau_{y} < \tau_1\}$$
$$\leq \varepsilon + \mathbf{P}_{y_1}\{\tau_{\lambda} < \tau_1\}.$$

Since $\mathbf{P}_{y_1}\{\tau_{\lambda} < \tau_1\}$ converges to 0 as $\lambda \to \infty$, there exists $\lambda_1 > \lambda_0$ such that $\mathbf{P}_{y_1}\{\tau_{\lambda} < \tau_1\} \leq \varepsilon$ for all $\lambda > \lambda_1$. Thus, we have proved that

$$\sup_{y \geq 1} \mathbf{P}_y\{\tau_{y} < \tau_1\} \leq 2\varepsilon \qquad \text{for all } \lambda > \lambda_1. \qquad \square$$

## 5. Excursion of $Z$.

5.1. *The local time at level* 1 *and the associated excursion measure.* According to Lemmas 3.1 and 3.2, the Markov process $Z$ is a Feller process in the filtration $\mathcal{F}$ for which 1 is regular for itself and instantaneous. It is therefore a "nice" Markov process in the sense of Chapter IV of [1] and we may consider a local time process $(L_t, t \geq 0)$ of $Z$ at level 1. Precisely, the local time process is such that:

- $(L_t, t \geq 0)$ is a continuous, $\mathcal{F}$-adapted process which increases on the closure of the set $\{t \geq 0, Z_t = 1\}$.
- For any stopping time $T$ such that $Z_T = 1$ a.s., the shifted process $(Z_{t+T}, L_{T+t} - L_T)_{t \geq 0}$ is independent of $\mathcal{F}_t$ and has the same law as $(Z_t, L_t)_{t \geq 0}$ under $\mathbf{P}_1$.

We can also consider the associated excursion measure $\mathbf{n}$ of the excursions of $Z$ away from 1 which we define in Chapter IV.4 of [1]. We denote by $(\epsilon_t, 0 \leq t \leq \zeta)$ a generic excursion with lifetime $\zeta$. Let also $L^{-1}$ stand for the right continuous inverse of $L$:

$$(5.1) \qquad L_t^{-1} \overset{\text{def}}{=} \inf\{s \geq 0, L_s > t\} \qquad \text{for all } t \geq 0.$$

Note that $L_t^{-1} < \infty$ for all $t$ since $Z$ is recurrent.

LEMMA 5.1. *Under $\mathbf{P}_1$, the process $L^{-1}$ is a subordinator whose Laplace exponent $\varphi$ defined by $\mathbf{E}_1[e^{-\lambda L_t^{-1}}] \stackrel{\text{def}}{=} e^{-t\varphi(\lambda)}$ has the form*

$$\varphi(\lambda) = \lambda \int_0^\infty e^{-\lambda r} \mathbf{n}\{\zeta > r\}\, dr.$$

*Moreover, there exist $c_{22}, c_{23} > 0$ such that $\mathbf{n}\{\zeta > r\} \leq c_{22} e^{-c_{23} r}$ for all $r \geq 1$, in particular, $\mathbf{n}[\zeta] < \infty$.*

PROOF. According to Theorem 8, page 114 of [1], $L^{-1}$ is a subordinator and its Laplace exponent $\varphi$ has the form

(5.2) $$\varphi(\lambda) = \lambda \mathbf{d} + \lambda \int_0^\infty e^{-\lambda r} \mathbf{n}\{\zeta > r\}\, dr.$$

Moreover, the drift coefficient $\mathbf{d}$ is such that $\mathbf{d}L(t) = \int_0^t \mathbf{1}_{\{Z_t = 1\}}\, dt$ $\mathbf{P}_1$-a.s. (cf. Corollary 6, page 112 of [1]). Thus, Lemma 3.3 implies that $\mathbf{d} = 0$. We now estimate the tail distribution of $\zeta$ under $\mathbf{n}$. Recall that $\tau_A(\epsilon)$ stands for the hitting time of the set $A$ for the excursion $\epsilon$:

$$\tau_A(\epsilon) \stackrel{\text{def}}{=} \inf\{t \in [0, \zeta], \epsilon_t \in A\} \qquad \text{(with the convention } \inf \varnothing = \infty\text{)}.$$

Since a generic excursion $\epsilon$ has no positive jumps, the Markov property yields, for $r > 1$,

$$\begin{aligned}
\mathbf{n}\{\zeta > r\} &\leq \mathbf{n}\{\tau_2(\epsilon) \leq 1, \zeta > r\} + \mathbf{n}\{\epsilon_1 \leq 2, \zeta > r\} \\
&\leq \mathbf{n}\{\tau_2(\epsilon) \leq 1, \zeta > 1\}\mathbf{P}_2\{\tau_1 > r - 1\} \\
&\quad + \mathbf{n}\{\epsilon_1 \leq 2, \zeta > 1\} \sup_{x \in (0,2)} \mathbf{P}_x\{\tau_1 > r - 1\} \\
&\leq 2\mathbf{n}\{\zeta > 1\} \sup_{x \in (0,2]} \mathbf{P}_x\{\tau_1 > r - 1\}.
\end{aligned}$$

Combining Lemma 3.4 and Proposition 4.3, we also have

$$\sup_{x \in (0,2]} \mathbf{P}_x\{\tau_1 > r - 1\} \leq \max(\mathbf{P}_0\{\tau_1 > r - 1\}, \mathbf{P}_2\{\tau_1 > r - 1\})$$

$$\leq c_{24} e^{-c_{25}(r-1)}.$$

This yields our estimate for $\mathbf{n}\{\zeta > r\}$. Finally, any excursion measure fulfills $\int_0^1 \mathbf{n}\{\zeta > r\}\, dr < \infty$, thus, $\mathbf{n}[\zeta] = \int_0^\infty \mathbf{n}\{\zeta > r\}\, dr < \infty$. $\quad\square$

LEMMA 5.2. *Let $f$ be a nonnegative measurable function. For all $\lambda > 0$, we have:*

(a)     $\mathbf{E}_1\left[\int_0^\infty e^{-\lambda t} f(Z_t)\, dt\right] = \dfrac{1}{\varphi(\lambda)} \mathbf{n}\left[\int_0^\zeta e^{-\lambda t} f(\epsilon_t)\, dt\right],$

(b)     $\mathbf{E}_1\left[\left(\int_0^\infty e^{-\lambda t} f(Z_t)\, dt\right)^2\right]$

$$= \frac{1}{\varphi(2\lambda)} \mathbf{n}\left[\left(\int_0^\zeta e^{-\lambda t} f(\epsilon_t)\, dt\right)^2\right]$$

$$+ \frac{2}{\varphi(\lambda)\varphi(2\lambda)} \mathbf{n}\left[\int_0^\zeta e^{-\lambda t} f(\epsilon_t)\, dt\right] \mathbf{n}\left[e^{-\lambda\zeta} \int_0^\zeta e^{-\lambda t} f(\epsilon_t)\, dt\right].$$

PROOF. Assertion (a) is a direct application of the compensation formula in excursion theory combined with the fact that the set $\{t \geq 0, Z_t = 1\}$ has 0 Lebesgue measure under $\mathbf{P}_1$ (Lemma 3.3). Compare with the example on page 120 of [1] for details.

We now prove (b). We use the notation $G_\lambda f(x) = \mathbf{E}_x[\int_0^\infty e^{-\lambda t} f(Z_t)\, dt]$. From a change of variable and with the help of the Markov property of $Z$,

$$\mathbf{E}_1\left[\left(\int_0^\infty e^{-\lambda t} f(Z_t)\, dt\right)^2\right] = 2\mathbf{E}_1\left[\int_0^\infty e^{-\lambda t} f(Z_t) \int_t^\infty e^{-\lambda s} f(Z_s)\, ds\, dt\right]$$

$$= 2\mathbf{E}_1\left[\int_0^\infty e^{-2\lambda t} f(Z_t) G_\lambda f(Z_t)\, dt\right].$$

Thus, using (a) with the function $x \mapsto f(x) G_\lambda f(x)$, we get that

$$(5.3) \quad \mathbf{E}_1\left[\left(\int_0^\infty e^{-\lambda t} f(Z_t)\, dt\right)^2\right] = \frac{2}{\varphi(2\lambda)} \mathbf{n}\left[\int_0^\zeta e^{-2\lambda t} f(\epsilon_t) G_\lambda f(\epsilon_t)\, dt\right].$$

We also have, with the help of the Markov property,

$$G_\lambda f(z) = \mathbf{E}_z\left[\int_0^{\tau_1} e^{-\lambda s} f(Z_s)\, ds\right] + \mathbf{E}_z\left[\int_{\tau_1}^\infty e^{-\lambda s} f(Z_s)\, ds\right]$$

$$= \mathbf{E}_z\left[\int_0^{\tau_1} e^{-\lambda s} f(Z_s)\, ds\right] + \mathbf{E}_z[e^{-\lambda\tau_1}] G_\lambda f(1).$$

Therefore, we may rewrite (5.3) as

$$(5.4) \quad \frac{2}{\varphi(2\lambda)} \mathbf{n}\left[\int_0^\zeta e^{-2\lambda t} f(\epsilon_t) \mathbb{E}_{\epsilon_t}\left[\int_0^{\tau_1} e^{-\lambda s} f(Z_s)\, ds\right] dt\right]$$

$$+ \frac{2}{\varphi(2\lambda)} \mathbf{n}\left[\int_0^\zeta e^{-2\lambda t} f(\epsilon_t) \mathbf{E}_{\epsilon_t}[e^{-\lambda\tau_1}]\, dt\right] G_\lambda f(1).$$

We deal with each term separately. Making use of the Markov property of the excursion $\epsilon$ at time $t$ under $\mathbf{n}(\cdot|\zeta > t)$ and with a change of variable, the first term of the last sum is equal to

$$(5.5) \quad \frac{2}{\varphi(2\lambda)} \mathbf{n} \left[ \int_0^\zeta e^{-2\lambda t} f(\epsilon_t) \int_t^\zeta e^{-\lambda(s-t)} f(\epsilon_s) \, ds \, dt \right]$$

$$= \frac{1}{\varphi(2\lambda)} \mathbf{n} \left[ \left( \int_0^\zeta e^{-\lambda t} f(\epsilon_t) \, dt \right)^2 \right].$$

Similarly, the second term of (5.4) may be rewritten

$$(5.6) \quad \frac{2G_\lambda f(1)}{\varphi(2\lambda)} \mathbf{n} \left[ \int_0^\zeta e^{-2\lambda t} f(\epsilon_t) e^{-\lambda(\zeta-t)} \, dt \right]$$

$$= \frac{2}{\varphi(\lambda)\varphi(2\lambda)} \mathbf{n} \left[ \int_0^\zeta e^{-\lambda t} f(\epsilon_t) \, dt \right] \mathbf{n} \left[ e^{-\lambda \zeta} \int_0^\zeta e^{-\lambda t} f(\epsilon_t) \, dt \right],$$

where we used (a) for the expression of $G_\lambda f(1)$ for the last equality. The combination of (5.3), (5.4), (5.5) and (5.6) yields (b). $\quad \square$

**COROLLARY 5.1.** *Let $g$ be a measurable, nonnegative function which is continuous almost everywhere with respect to the Lebesgue measure. Then*

$$\mathbf{n} \left[ \int_0^\zeta g(\epsilon_t) \, dt \right] = \mathbf{n}[\zeta] \mathbf{E}[g(Z_\infty)].$$

PROOF. In view of the monotone convergence theorem, we may assume that $g$ is bounded. First, using (a) of the previous lemma with the function $f = 1$,

$$(5.7) \quad \frac{\varphi(\lambda)}{\lambda} = \mathbf{n} \left[ \int_0^\zeta e^{-\lambda t} \, dt \right] \underset{\lambda \to 0+}{\longrightarrow} \mathbf{n}[\zeta].$$

Thus, using again (a) of Lemma 5.2 but now with the function $g$, and with the help of the monotone convergence theorem, we find

$$\mathbf{n} \left[ \int_0^\zeta g(\epsilon_t) \, dt \right] = \lim_{\lambda \to 0+} \varphi(\lambda) \mathbf{E}_1 \left[ \int_0^\infty e^{-\lambda t} g(Z_t) \, dt \right]$$

$$= \mathbf{n}[\zeta] \lim_{\lambda \to 0+} \mathbf{E}_1 \left[ \lambda \int_0^\infty e^{-\lambda t} g(Z_t) \, dt \right].$$

By a change of variable and using Fubini's theorem, we also have

$$\mathbf{E}_1 \left[ \lambda \int_0^\infty e^{-\lambda t} g(Z_t) \, dt \right] = \int_0^\infty \mathbf{E}_1 [g(Z_{y/\lambda})] e^{-y} \, dy.$$

For any $y > 0$, $Z_{y/\lambda}$ converges in law toward $Z_\infty$ as $\lambda \to 0+$. Moreover, according to Proposition 3.1, $Z_\infty$ has a continuous density with respect to the Lebesgue measure and $g$ is continuous almost everywhere, hence, $\lim_{\lambda \to 0+} \mathbf{E}_1[g(Z_{y/\lambda})] = \mathbb{E}[g(Z_\infty)]$. Making use of the dominated convergence theorem, we conclude that

$$\mathbf{E}_1\left[\lambda \int_0^\infty e^{-\lambda t} g(Z_t)\, dt\right] \underset{\lambda \to 0+}{\longrightarrow} \int_0^\infty \mathbf{E}[g(Z_\infty)] e^{-y}\, dy = \mathbf{E}[g(Z_\infty)]. \qquad \square$$

COROLLARY 5.2. *Recall that* $\mathbf{m} \overset{\text{def}}{=} \frac{-2}{\Phi(1)}$. *When* $\kappa > 1$ *[i.e., when* $\Phi(1) < 0$*], we have*

$$\mathbf{n}\left[\int_0^\zeta \epsilon_t\, dt\right] = \mathbf{n}[\zeta]\mathbf{m}.$$

PROOF. Corollary 5.1 yields $\mathbf{n}[\int_0^\zeta \epsilon_t\, dt] = \mathbf{n}[\zeta]\mathbf{E}[Z_\infty]$. According to Proposition 3.1, $Z_\infty$ has the same law as $U(1) \int_0^\infty e^{\mathbb{V}_s}\, ds$ under $\mathbf{P}_0$. Moreover, $U(1)$ under $\mathbf{P}_0$ has an exponential distribution with mean 2 and is independent of $\mathbb{V}$, hence,

$$\mathbf{E}[Z_\infty] = 2\int_0^\infty \mathbf{E}[e^{\mathbb{V}_s}]\, ds = 2\int_0^\infty e^{t\Phi(1)}\, ds = -\frac{2}{\Phi(1)}. \qquad \square$$

5.2. *Maximum of an excursion.* The aim of this subsection is to study the distribution of the supremum of an excursion. Our main result is contained in the following proposition.

PROPOSITION 5.1. *We have*

$$\mathbf{n}\{\tau_z(\epsilon) < \infty\} \underset{z \to \infty}{\sim} \mathbf{n}[\zeta]\frac{2^\kappa \Gamma(\kappa)\kappa^2 \mathbf{K}}{z^\kappa}.$$

Of course, this estimate may be rewritten

$$\mathbf{n}\left\{\sup_{[0,\zeta]} \epsilon > z\right\} \underset{z \to \infty}{\sim} \mathbf{n}[\zeta]\frac{2^\kappa \Gamma(\kappa)\kappa^2 \mathbf{K}}{z^\kappa}.$$

The proof relies on two lemmas.

LEMMA 5.3. *We have*

$$\mathbf{E}\left[\int_0^\infty \mathbf{1}_{\{\mathbb{V}_t > 0\}}\, dt\right] = \mathbf{E}\left[\int_0^\infty \mathbf{1}_{\{\mathbb{V}_t \geq 0\}}\, dt\right] = \frac{1}{\kappa\Phi'(\kappa)}.$$

LEMMA 5.4. *We have*

$$\lim_{z \to \infty} \mathbf{E}_z\left[\int_0^{\tau_1} \mathbf{1}_{\{Z_t \geq z\}}\, dt\right] = \frac{1}{\kappa\Phi'(\kappa)}.$$

Let us for the time being admit the lemmas and give the proof of the proposition.

PROOF OF PROPOSITION 5.1. Since a generic excursion $\epsilon$ under $\mathbf{n}$ has no positive jumps, the Markov property yields

$$
(5.8) \quad \begin{aligned}
\mathbf{n}\left[\int_0^\zeta \mathbf{1}_{\{\epsilon_s > z\}} \, ds\right] &= \mathbf{n}\left[\mathbf{1}_{\{\tau_z(\epsilon) < \infty\}} \int_{\tau_z(\epsilon)}^\zeta \mathbf{1}_{\{\epsilon_s > z\}} \, ds\right] \\
&= \mathbf{n}\{\tau_z(\epsilon) < \infty\} \mathbb{E}_z\left[\int_0^{\tau_1} \mathbf{1}_{\{Z_s > z\}} \, ds\right].
\end{aligned}
$$

On the one hand, from Corollary 5.1 and Proposition 3.1,

$$
(5.9) \quad \mathbf{n}\left[\int_0^\zeta \mathbf{1}_{\{\epsilon_s > z\}} \, ds\right] = \mathbf{n}[\zeta]\mathbf{P}\{Z_\infty > z\} \underset{z \to \infty}{\sim} \mathbf{n}[\zeta]\frac{2^\kappa \Gamma(\kappa+1)\mathbf{K}}{\Phi'(\kappa)z^\kappa}.
$$

On the other hand, according to Lemma 5.4,

$$
(5.10) \quad \mathbb{E}_z\left[\int_0^{\tau_1} \mathbf{1}_{\{Z_s > z\}} \, ds\right] \underset{z \to \infty}{\longrightarrow} \frac{1}{\kappa\Phi'(\kappa)}.
$$

The proposition follows from the combination of (5.8), (5.9) and (5.10).  □

PROOF OF LEMMA 5.3. Since $\mathbb{V}$ has no positive jumps, it is not a compound Poisson process, therefore, Proposition 15, page 30 of [1] states that the resolvent measures of $\mathbb{V}$ are diffuse, that is, $\mathbf{E}[\int_0^\infty \mathbf{1}_{\{\mathbb{V}_t = 0\}} \, dt] = 0$. Thus,

$$
\mathbf{E}\left[\int_0^\infty \mathbf{1}_{\{\mathbb{V}_t \geq 0\}} \, dt\right] = \mathbf{E}\left[\int_0^\infty \mathbf{1}_{\{\mathbb{V}_t > 0\}} \, dt\right].
$$

Let $\Psi : [0, \infty) \mapsto [\kappa, \infty)$ denote the right inverse of the Laplace exponent $\Phi$ such that $\Phi \circ \Psi(\lambda) = \lambda$ for all $\lambda \geq 0$ [in particular, $\Psi(0) = \kappa$]. Then, Exercise 1 on page 212 of [1] which is an easy consequence of Corollary 3, page 190 of [1] states that

$$
\mathbf{E}\left[\int_0^\infty e^{-\lambda t} \mathbf{1}_{\{\mathbb{V}_t \geq 0\}} \, dt\right] = \frac{\Psi'(\lambda)}{\Psi(\lambda)} \qquad \text{for all } \lambda > 0.
$$

Taking the limit as $\lambda \to 0$, we conclude that

$$
\mathbf{E}\left[\int_0^\infty \mathbf{1}_{\{\mathbb{V}_t \geq 0\}} \, dt\right] = \frac{\Psi'(0)}{\Psi(0)} = \frac{1}{\kappa\Phi'(\kappa)}. \qquad\qquad □
$$

PROOF OF LEMMA 5.4. Assume that $z > 1$ and let $\varepsilon > 0$. Note that for $1 < b < z$, we have $\tau_{[0,z/b]} \leq \tau_1$ $\mathbf{P}_z$-a.s. Thus, on the one hand,

$$
(5.11) \quad \mathbf{E}_z\left[\int_0^{\tau_{[0,z/b]}} \mathbf{1}_{\{Z_t \geq z\}} \, dt\right] \leq \mathbf{E}_z\left[\int_0^{\tau_1} \mathbf{1}_{\{Z_t \geq z\}} \, dt\right].
$$

On the other hand, making use of the Markov property of $Z$ and with the help of Lemma 3.4,

$$(5.12)\quad
\begin{aligned}
&\mathbf{E}_z\left[\int_0^{\tau_1}\mathbf{1}_{\{Z_t\geq z\}}\,dt\right]\\
&= \mathbf{E}_z\left[\int_0^{\tau_{[0,z/b]}}\mathbf{1}_{\{Z_t\geq z\}}\,dt\right] + \mathbf{E}_z\left[\int_{\tau_{[0,z/b]}}^{\tau_1}\mathbf{1}_{\{Z_t\geq z\}}\,dt\right]\\
&\leq \mathbf{E}_z\left[\int_0^{\tau_{[0,z/b]}}\mathbf{1}_{\{Z_t\geq z\}}\,dt\right] + \mathbf{E}_{z/b}\left[\int_0^{\tau_1}\mathbf{1}_{\{Z_t\geq z\}}\,dt\right]\\
&= \mathbf{E}_z\left[\int_0^{\tau_{[0,z/b]}}\mathbf{1}_{\{Z_t\geq z\}}\,dt\right] + \mathbf{P}_{z/b}\{\tau_z<\tau_1\}\mathbf{E}_z\left[\int_0^{\tau_1}\mathbf{1}_{\{Z_t\geq z\}}\,dt\right].
\end{aligned}
$$

According to Proposition 4.4, there exists $b_1>1$ such that, for $b>b_1$, $\sup_{z\geq b}\mathbf{P}_{z/b}\{\tau_z<\tau_1\}\leq\varepsilon$. Therefore, combining (5.11) and (5.12), for all $z>b>b_1$,

$$
\begin{aligned}
\mathbf{E}_z\left[\int_0^{\tau_{[0,z/b]}}\mathbf{1}_{\{Z_t\geq z\}}\,dt\right] &\leq \mathbf{E}_z\left[\int_0^{\tau_1}\mathbf{1}_{\{Z_t\geq z\}}\,dt\right]\\
&\leq \frac{1}{1-\varepsilon}\mathbf{E}_z\left[\int_0^{\tau_{[0,z/b]}}\mathbf{1}_{\{Z_t\geq z\}}\,dt\right].
\end{aligned}
$$

Thus, we just need to prove that we may find $b_2>b_1$ and $z_0>0$ such that

$$(5.13)\quad
\begin{aligned}
\frac{1}{\kappa\Phi'(\kappa)}-\varepsilon &\leq \mathbf{E}_z\left[\int_0^{\tau_{[0,z/b_2]}}\mathbf{1}_{\{Z_t\geq z\}}\,dt\right]\\
&\leq \frac{1}{\kappa\Phi'(\kappa)}+\varepsilon \qquad\text{for all }z\geq z_0.
\end{aligned}
$$

Recall from Lemma 3.1 that $Z$ under $\mathbf{P}_z$ has the same law as the process $(ze^{\mathbb{V}_t}U(a(t)/z),t\geq 0)$ under $\mathbf{P}_1$. Thus,

$$(5.14)\quad
\begin{aligned}
&\mathbf{P}_z\{Z_t\geq z,\tau_{[0,z/b]}\geq t\}\\
&= \mathbf{P}_1\left\{e^{\mathbb{V}_t}U\left(\frac{a(t)}{z}\right)\geq 1,\forall s\in[0,t)\ e^{\mathbb{V}_s}U\left(\frac{a(s)}{z}\right)>\frac{1}{b}\right\}.
\end{aligned}
$$

Since $U$ is continuous at 0 and starting from 1 under $\mathbb{P}_1$, we also have

$$(5.15)\quad
\sup_{0\leq s\leq t}\left|U\left(\frac{a(s)}{z}\right)-1\right|\xrightarrow[z\to\infty]{\mathbf{P}_1\text{-a.s.}}0 \qquad\text{for all }t\geq 0.
$$

Combining (5.14) and (5.15), we get that, for all fixed $t\geq 0$,

$$
\begin{aligned}
\liminf_{z\to\infty}\mathbf{P}_z\{Z_t\geq z,\tau_{[0,z/b]}\geq t\} &\geq \mathbf{P}\left\{e^{\mathbb{V}_t}>1,\forall s\in[0,t)\ e^{\mathbb{V}_s}>\frac{1}{b}\right\}\\
&= \mathbf{P}\{\mathbb{V}_t>0,\tau_{(-\infty,-\log b]}(\mathbb{V})\geq t\}.
\end{aligned}
$$

Thus, by inversion of the sum and from Fatou's lemma,

$$
\liminf_{z \to \infty} \mathbf{E}_z \left[ \int_0^{\tau_{[0,z/b]}} \mathbf{1}_{\{Z_t \geq z\}} \, dt \right] = \liminf_{z \to \infty} \int_0^\infty \mathbf{P}_z \{ Z_t \geq z, \tau_{[0,z/b]} \geq t \} \, dt
$$
$$
\geq \int_0^\infty \liminf_{z \to \infty} \mathbf{P}_z \{ Z_t \geq z, \tau_{[0,z/b]} \geq t \} \, dt
$$
$$
\geq \int_0^\infty \mathbf{P} \{ \mathbb{V}_t > 0, \tau_{(-\infty, -\log b]}(\mathbb{V}) \geq t \} \, dt
$$
$$
= \mathbf{E} \left[ \int_0^{\tau_{(-\infty, -\log b]}(\mathbb{V})} \mathbf{1}_{\{\mathbb{V}_t > 0\}} \, dt \right].
$$

By use of the monotone convergence theorem, we also have

$$
\lim_{b \to \infty} \mathbf{E} \left[ \int_0^{\tau_{(-\infty, -\log b]}(\mathbb{V})} \mathbf{1}_{\{\mathbb{V}_t > 0\}} \, dt \right] = \mathbf{E} \left[ \int_0^\infty \mathbf{1}_{\{\mathbb{V}_t > 0\}} \, dt \right] = \frac{1}{\kappa \Phi'(\kappa)},
$$

where we used Lemma 5.3 for the last equality. We may therefore find $b_2 > b_1$ such that, for all $z$ large enough,

$$
\mathbf{E}_z \left[ \int_0^{\tau_{[0,z/b_2]}} \mathbf{1}_{\{Z_t \geq z\}} \, dt \right] \geq \frac{1}{\kappa \Phi'(\kappa)} - \varepsilon.
$$

We still have to prove the upper bound in (5.13). Keeping in mind (5.15), we notice that, for all fixed $t \geq 0$,

$$
\limsup_{z \to \infty} \mathbf{P}_z \{ Z_t \geq z, \tau_{[0,z/b_2]} \geq t \} \leq \limsup_{z \to \infty} \mathbf{P}_z \{ Z_t \geq z \} \leq \mathbf{P} \{ \mathbb{V}_t \geq 0 \}.
$$

Moreover, Proposition 4.2 states that there exist $c_{26,b_2}, c_{27,b_2} > 0$ such that, for all $z$ large enough and all $t \geq 0$,

$$
\mathbf{P}_z \{ Z_t \geq z, \tau_{[0,z/b_2]} \geq t \} \leq \mathbf{P}_z \{ \tau_{[0,z/b_2]} > t \} \leq c_{26,b_2} e^{-c_{27,b_2} t}.
$$

This domination result enables us to use Fatou's lemma for the lim sup. Thus, just as for the lim inf, we now find

$$
\limsup_{z \to \infty} \mathbf{E}_z \left[ \int_0^{\tau_{[0,z/b_2]}} \mathbf{1}_{\{Z_t \geq z\}} \, dt \right] \leq \int_0^\infty \limsup_{z \to \infty} \mathbf{P}_z \{ Z_t \geq z, \tau_{[0,z/b_2]} \geq t \} \, dt
$$
$$
\leq \int_0^\infty \mathbf{P} \{ \mathbb{V}_t \geq 0 \} \, dt
$$
$$
= \mathbf{E} \left[ \int_0^\infty \mathbf{1}_{\{\mathbb{V}_t \geq 0\}} \, dt \right] = \frac{1}{\kappa \Phi'(\kappa)}.
$$

This completes the proof of the lemma. $\square$

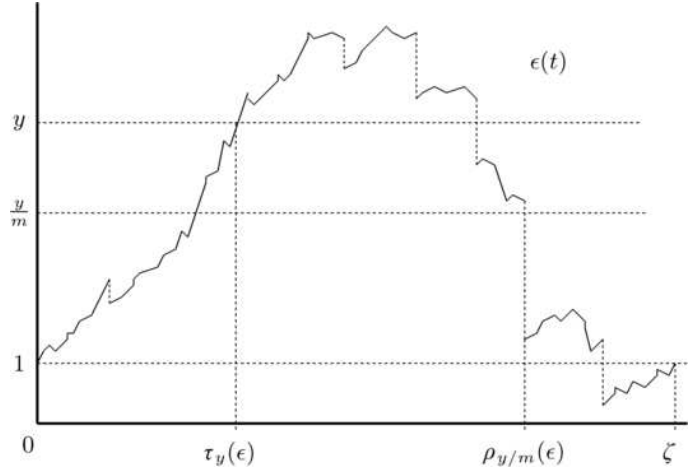

Fig. 3.    *An excursion $\epsilon$.*

5.3. *Integral of an excursion.* We now estimate the tail distribution of the area of an excursion. The next proposition is the key to the proof of our main theorem.

PROPOSITION 5.2.    *We have*

$$\mathbf{n}\left\{\int_0^\zeta \epsilon_s\,ds > x\right\} \underset{x\to\infty}{\sim} \mathbf{n}[\zeta]\frac{2^\kappa \Gamma(\kappa)\kappa^2 \mathbf{K}^2}{\Phi'(\kappa)x^\kappa}.$$

In the rest of this subsection, we assume $x$ to be a large number and we will use the notation

$$(5.16) \qquad\qquad m \overset{\text{def}}{=} \log^3 x,$$

$$(5.17) \qquad\qquad y \overset{\text{def}}{=} \frac{x}{m} = \frac{x}{\log^3 x}.$$

The idea of the proof of the proposition is to decompose the integral of an excursion $\epsilon$ such that $\tau_y(\varepsilon) < \infty$ in the form (see Figure 3)

$$(5.18) \qquad \int_0^\zeta \epsilon_s\,ds = \int_0^{\tau_y(\epsilon)} \epsilon_s\,ds + \int_{\tau_y(\varepsilon)}^{\rho_{y/m}(\epsilon)} \epsilon_s\,ds + \int_{\rho_{y/m}(\epsilon)}^\zeta \epsilon_s\,ds,$$

where $\rho_{y/m} = \inf\{t > \tau_y(\epsilon), \epsilon_t \le y/m\}$. We will show that the contributions of the first and last term on the r.h.s. of (5.18) are negligible. As for the second term, we will show that its distribution is well approximated by the distribution of the random variable $y\int_0^\infty e^{\mathbb{V}_t}\,dt$. This will give

$$\mathbf{n}\left\{\int_0^\zeta \epsilon_s\,ds > x\right\} \approx \mathbf{n}\{\tau_y(\epsilon) < \infty\}\mathbf{P}\left\{y\int_0^\infty e^{\mathbb{V}_t} > x\right\}$$

and the proposition will follow from the estimates obtained in the previous sections. We start with a lemma:

LEMMA 5.5. *Recall notation (5.16) and (5.17). We have*

$$\mathbf{P}_y\left\{\int_0^{\tau_{[0,y/m]}} Z_s\, ds > x\right\} \underset{x\to\infty}{\sim} \frac{\mathbf{K}}{\Phi'(\kappa)}\left(\frac{y}{x}\right)^\kappa.$$

PROOF. Let $(\widetilde{Z}_t, t \geq 0)$ denote the process

$$\widetilde{Z}_t = ye^{\mathbb{V}_t}U\left(\frac{a(t)}{y}\right).$$

We have already proved in Lemma 3.1 that $\widetilde{Z}$ under $\mathbf{P}_1$ has the same law as $Z$ under $\mathbf{P}_y$. Let $\widetilde{\tau}_A$ denote the hitting time of the set $A$ for the process $\widetilde{Z}$. We must prove that

$$\mathbf{P}_1\left\{\int_0^{\widetilde{\tau}_{[0,y/m]}} e^{\mathbb{V}_t}U\left(\frac{a(t)}{y}\right)dt > \frac{x}{y}\right\} \underset{x\to\infty}{\sim} \frac{\mathbf{K}}{\Phi'(\kappa)}\left(\frac{y}{x}\right)^\kappa.$$

We define

$$\gamma \overset{\text{def}}{=} \inf\{t \geq 0, \mathbb{V}_t < -\log(2m)\},$$

$$\gamma' \overset{\text{def}}{=} \inf\{t \geq 0, \mathbb{V}_t < -\log(m/2)\},$$

and for $0 < \varepsilon < \frac{1}{2}$, set

$$\mathcal{E} \overset{\text{def}}{=} \left\{|U(z) - 1| \leq \varepsilon \text{ for all } 0 \leq z \leq \frac{2m\gamma}{y}\right\}.$$

Let us first notice that for all $0 \leq t \leq \gamma$, we have $a(t) = \int_0^t e^{-\mathbb{V}_s}\, ds \leq 2m\gamma$ and $e^{\mathbb{V}_\gamma} \leq \frac{1}{2m}$. Thus, on $\mathcal{E}$, we have

$$(5.19) \qquad \widetilde{Z}_\gamma = ye^{\mathbb{V}_\gamma}U\left(\frac{a(\gamma)}{y}\right) < \frac{y}{2m}(1+\varepsilon) < \frac{y}{m}.$$

We also have $e^{\mathbb{V}_t} \geq \frac{2}{m}$ for all $t < \gamma' \leq \gamma$. Thus, on $\mathcal{E}$,

$$(5.20) \qquad \widetilde{Z}_t = ye^{\mathbb{V}_t}U\left(\frac{a(t)}{y}\right) \geq \frac{2y}{m}(1-\varepsilon) > \frac{y}{m} \qquad \text{for all } t < \gamma'.$$

Combining (5.19) and (5.20), we deduce that

$$\mathcal{E} \subset \{\gamma' \leq \widetilde{\tau}_{[0,y/m]} \leq \gamma\}.$$

Let us for the time being admit that

$$(5.21) \qquad \lim_{x\to\infty}\left(\frac{x}{y}\right)^\kappa \mathbf{P}_1\{\mathcal{E}^c\} = 0.$$

We now write

$$\mathbf{P}_1\left\{\int_0^{\widetilde{\tau}_{[0,y/m]}} e^{\mathbb{V}_t} U\left(\frac{a(t)}{y}\right) dt > \frac{x}{y}\right\}$$

$$\leq \mathbf{P}_1\left\{\int_0^{\widetilde{\tau}_{[0,y/m]}} e^{\mathbb{V}_t} U\left(\frac{a(t)}{y}\right) dt > \frac{x}{y}, \mathcal{E}\right\} + \mathbf{P}_1\{\mathcal{E}^c\}$$

$$\leq \mathbf{P}_1\left\{\int_0^{\widetilde{\tau}_{[0,y/m]}} e^{\mathbb{V}_t}(1+\varepsilon)\, dt > \frac{x}{y}, \mathcal{E}\right\} + \mathbf{P}_1\{\mathcal{E}^c\}$$

$$\leq \mathbf{P}\left\{\int_0^{\infty} e^{\mathbb{V}_t}\, dt > \frac{x}{(1+\varepsilon)y}\right\} + \mathbf{P}_1\{\mathcal{E}^c\}.$$

We have already checked in the proof of Proposition 3.1 that

$$\mathbf{P}\left\{\int_0^{\infty} e^{\mathbb{V}_t}\, dt > \frac{x}{(1+\varepsilon)y}\right\} \underset{x\to\infty}{\sim} \frac{\mathbf{K}}{\Phi'(\kappa)}\left(\frac{(1+\varepsilon)y}{x}\right)^{\kappa}.$$

Therefore,

$$\limsup_{x\to\infty}\left(\frac{x}{y}\right)^{\kappa}\mathbf{P}_1\left\{\int_0^{\widetilde{\tau}_{[0,y/m]}} e^{\mathbb{V}_t} U\left(\frac{a(t)}{y}\right) dt > \frac{x}{y}\right\} \leq \frac{\mathbf{K}(1+\varepsilon)^{\kappa}}{\Phi'(\kappa)}.$$

We now prove the lim inf. Since $\gamma' \leq \widetilde{\tau}_{[0,y/m]}$ on $\mathcal{E}$,

$$\mathbf{P}_1\left\{\int_0^{\widetilde{\tau}_{[0,y/m]}} e^{\mathbb{V}_t} U\left(\frac{a(t)}{y}\right) dt > \frac{x}{y}\right\}$$

$$\geq \mathbf{P}_1\left\{\int_0^{\widetilde{\tau}_{[0,y/m]}} e^{\mathbb{V}_t} U\left(\frac{a(t)}{y}\right) dt > \frac{x}{y}, \mathcal{E}\right\} - \mathbf{P}_1\{\mathcal{E}^c\}$$

$$\geq \mathbf{P}_1\left\{\int_0^{\gamma'} e^{\mathbb{V}_t}(1-\varepsilon)\, dt > \frac{x}{y}, \mathcal{E}\right\} - \mathbf{P}_1\{\mathcal{E}^c\}$$

$$\geq \mathbf{P}\left\{\int_0^{\gamma'} e^{\mathbb{V}_t}\, dt > \frac{x}{(1-\varepsilon)y}\right\} - 2\mathbf{P}_1\{\mathcal{E}^c\}.$$

Since $\mathbb{V}_{\gamma'} < -\log(m/2)$, it is easy to check with the help of the Markov property of $\mathbb{V}$ that

$$\mathbf{P}\left\{\int_0^{\gamma'} e^{\mathbb{V}_t}\, dt > \frac{x}{(1-\varepsilon)y}\right\} \underset{x\to\infty}{\sim} \mathbf{P}\left\{\int_0^{\infty} e^{\mathbb{V}_t}\, dt > \frac{x}{(1-\varepsilon)y}\right\}$$

$$\underset{x\to\infty}{\sim} \frac{\mathbf{K}}{\Phi'(\kappa)}\left(\frac{(1-\varepsilon)y}{x}\right)^{\kappa},$$

so we obtain the lower bound

$$\liminf_{x\to\infty}\left(\frac{x}{y}\right)^{\kappa}\mathbf{P}_1\left\{\int_0^{\widetilde{\tau}_{[0,y/m]}} e^{\mathbb{V}_t} U\left(\frac{a(t)}{y}\right) dt > \frac{x}{y}\right\} \geq \frac{\mathbf{K}(1-\varepsilon)^{\kappa}}{\Phi'(\kappa)}.$$

It remains to prove (5.21). To this end, notice that

$$\mathbf{P}_1\{\mathcal{E}^c\} \leq \mathbf{P}\left\{\frac{2m\gamma}{y} \geq \frac{m^2}{y}\right\} + \mathbf{P}_1\left\{\sup_{z\in[0,m^2/y]}|U(z)-1| > \varepsilon\right\}$$

$$\leq \mathbf{P}\{\mathbb{V}_{m/2} \geq -\log(2m)\} + \mathbf{P}_1\left\{\sup_{z\in[0,m^2/y]}|U(z)-1| > \varepsilon\right\}.$$

Recall that $\Phi(\kappa/2) < 0$. Thus, on the one hand,

$$\mathbf{P}\{\mathbb{V}_{m/2} \geq -\log(2m)\} \leq (2m)^{\kappa/2}\mathbf{E}[e^{\kappa/2\mathbb{V}_{m/2}}]$$

$$= (2m)^{\kappa/2}e^{m/2\Phi(\kappa/2)} = o\left(\left(\frac{y}{x}\right)^\kappa\right).$$

On the other hand, $U$ under $\mathbf{P}_1$ is a squared Bessel process of dimension 2 starting from 1. Thus, it has the same law as $B^2 + \widetilde{B}^2 + 2B + 1$, where $B$ and $\widetilde{B}$ are two independent Brownian motions. Hence,

$$\mathbf{P}_1\left\{\sup_{[0,m^2/y]}|U-1| > \varepsilon\right\} \leq 2\mathbf{P}\left\{\sup_{[0,m^2/y]}|B|^2 > \frac{\varepsilon}{4}\right\} + \mathbf{P}\left\{\sup_{[0,m^2/y]}|B| > \frac{\varepsilon}{4}\right\}$$

$$\leq 3\mathbf{P}\left\{\sup_{[0,m^2/y]}|B| > \frac{\varepsilon}{4}\right\}.$$

Finally, from the exact distribution of $\sup_{[0,1]}|B|$ and the usual estimate on Gaussian tails,

$$\mathbf{P}_1\left\{\sup_{[0,t]}|B| > \right\} \leq \frac{2\sqrt{t}}{a}e^{-a^2/(2t)} \qquad \text{for all } a, t > 0.$$

Therefore,

$$\mathbf{P}_1\left\{\sup_{z\in[0,m^2/y]}|U(z)-1| > \varepsilon\right\} \leq \frac{24m}{\varepsilon\sqrt{y}}\exp\left(-\frac{\varepsilon^2 y}{32m^2}\right) = o\left(\left(\frac{y}{x}\right)^\kappa\right).$$

This completes the proof of the lemma □

PROOF OF PROPOSITION 5.2. We first deal with the lim inf; we have

$$\mathbf{n}\left\{\int_0^\zeta \epsilon_s\,ds > x\right\} \geq \mathbf{n}\left\{\tau_y(\epsilon) < \infty, \int_{\tau_y(\epsilon)}^{\tau_{[0,y/m]}(\epsilon)} \epsilon_s\,ds > x\right\}.$$

Using the Markov property and the fact that the excursion $\epsilon$ does not possess positive jumps, the r.h.s. of this inequality is equal to

$$\mathbf{n}\{\tau_y(\epsilon) < \infty\}\mathbf{P}_y\left\{\int_0^{\tau_{[0,y/m]}} Z_s\,ds > x\right\} \underset{x\to\infty}{\sim} \mathbf{n}[\zeta]\frac{2^\kappa\Gamma(\kappa)\kappa^2\mathbf{K}^2}{\Phi'(\kappa)x^\kappa},$$

where we used Lemma 5.5 and Proposition 5.1 for the equivalence. Therefore,

$$\liminf_{x \to \infty} x^\kappa \mathbf{n} \left\{ \int_0^\zeta \epsilon_s \, ds > x \right\} \geq \frac{2^\kappa \Gamma(\kappa) \kappa^2 \mathbf{K}^2}{\Phi'(\kappa)}.$$

We now prove the upper bound. Let $\varepsilon > 0$. We simply need to show that

$$\limsup_{x \to \infty} x^\kappa \mathbf{n} \left\{ \int_0^\zeta \epsilon_s \, ds > (1 + 2\varepsilon)x \right\} \leq \frac{2^\kappa \Gamma(\kappa) \kappa^2 \mathbf{K}^2}{\Phi'(\kappa)}.$$

According to Lemma 5.1, we have $\mathbf{n}\{\zeta \geq \log^2 x\} = o(x^{-\kappa})$, thus

$$\mathbf{n} \left\{ \int_0^\zeta \epsilon_s \, ds > (1 + 2\varepsilon)x \right\} = \mathbf{n} \left\{ \zeta < \log^2 x, \int_0^\zeta \epsilon_s \, ds > (1 + 2\varepsilon)x \right\} + o(x^{-\kappa}).$$

We also note that $\int_0^\zeta \epsilon_s \, ds \leq \zeta \sup_{s \in [0, \zeta]} \epsilon_s$. Since $y = x / \log^3 x$, we deduce that, for all $x$ large enough,

$$\left\{ \zeta < \log^2 x, \int_0^\zeta \epsilon_s \, ds > (1 + 2\varepsilon)x \right\}$$

$$= \left\{ \tau_y(\varepsilon) < \zeta < \log^2 x, \int_0^\zeta \epsilon_s \, ds > (1 + 2\varepsilon)x, \int_0^{\tau_y(\epsilon)} \epsilon_s \, ds < \varepsilon x \right\}$$

$$\subset \left\{ \tau_y(\varepsilon) < \zeta < \log^2 x, \int_{\tau_y(\epsilon)}^\zeta \epsilon_s \, ds > (1 + \varepsilon)x \right\}.$$

Thus, making use of the Markov property of $\epsilon$ for the stopping time $\tau_y(\epsilon)$,

$$\mathbf{n} \left\{ \int_0^\zeta \epsilon_s \, ds > (1 + 2\varepsilon)x \right\}$$

$$\leq \mathbf{n} \left\{ \tau_y(\varepsilon) < \zeta < \log^2 x, \int_{\tau_y(\epsilon)}^\zeta \epsilon_s \, ds > (1 + \varepsilon)x \right\} + o(x^{-\kappa})$$

$$\leq \mathbf{n} \{ \tau_y(\varepsilon) < \infty \} \mathbf{P}_y \left\{ \int_0^{\tau_1} Z_s \, ds > (1 + \varepsilon)x, \tau_1 < \log^2 x \right\} + o(x^\kappa).$$

In view of Proposition 5.1, it remains to prove that

$$\limsup_{x \to \infty} \left( \frac{x}{y} \right)^\kappa \mathbf{P}_y \left\{ \int_0^{\tau_1} Z_s \, ds > (1 + \varepsilon)x, \tau_1 < \log^2 x \right\} \leq \frac{\mathbf{K}}{\Phi'(\kappa)}.$$

We have

$$\mathbf{P}_y \left\{ \int_0^{\tau_1} Z_s \, ds > (1 + \varepsilon)x, \tau_1 < \log^2 x \right\}$$

$$\leq \mathbf{P}_y \left\{ \int_0^{\tau_{[0, y/m]}} Z_s \, ds > x \right\} + \mathbf{P}_y \left\{ \int_{\tau_{[0, y/m]}}^{\tau_1} Z_s \, ds > \varepsilon x, \tau_1 < \log^2 x \right\}.$$

On the one hand, according to Lemma 5.5,

$$\lim_{x\to\infty}\left(\frac{x}{y}\right)^{\kappa}\mathbf{P}_y\left\{\int_0^{\tau_{[0,y/m]}} Z_s\,ds > x\right\} = \frac{\mathbf{K}}{\Phi'(\kappa)}.$$

On the other hand,

$$\begin{aligned}
\mathbf{P}_y&\left\{\int_{\tau_{[0,y/m]}}^{\tau_1} Z_s\,ds > \varepsilon x, \tau_1 < \log^2 x\right\}\\
(5.22)\qquad &\le \mathbf{P}_y\left\{\sup_{s\in[\tau_{[0,y/m]},\tau_1]} Z_s > \frac{\varepsilon x}{\log^2 x}\right\}\\
&\le \mathbf{P}_{y/m}\{\tau_{\varepsilon x/\log^2 x} < \tau_1\},
\end{aligned}$$

where we used the Markov property of $Z$ for the stopping time $\tau_{[0,y/m]}$ combined with Lemma 3.4 and the absence of positive jumps for the last inequality. Since $\frac{y}{m} < \frac{x}{\log^2 x}$, we also notice that

$$\begin{aligned}
\mathbf{n}\{\tau_{\varepsilon x/\log^2 x}(\epsilon) < \infty\} &= \mathbf{n}\{\tau_{y/m}(\epsilon) < \tau_{\varepsilon x/\log^2 x}(\epsilon) < \infty\}\\
&= \mathbf{n}\{\tau_{y/m}(\epsilon) < \infty\}\mathbf{P}_{y/m}\{\tau_{\varepsilon x/\log^2 x} < \tau_1\}.
\end{aligned}$$

Therefore, (5.22) is also equal to

$$\frac{\mathbf{n}\{\tau_{\varepsilon x/\log^2 x}(\epsilon) < \infty\}}{\mathbf{n}\{\tau_{y/m}(\epsilon) < \infty\}} \underset{x\to\infty}{\sim} \left(\frac{y\log^2 x}{\varepsilon x m}\right)^{\kappa} = o\left(\left(\frac{y}{x}\right)^{\kappa}\right),$$

where we used Proposition 5.1 for the equivalence. This concludes the proof of the proposition. □

## 6. The second moment.

Recall that $\mathbf{m} = -2/\Phi(1)$. The aim of this section is to calculate the quantity $\mathbf{n}[(\int_0^\zeta (\epsilon_t - \mathbf{m})\,dt)^2]$ when $\kappa > 2$ in terms of the Laplace exponent $\Phi$ of $\mathbb{V}$. We start with the following:

LEMMA 6.1. *When $\kappa > 2$, for all $t, z \ge 0$,*

(a) $\quad\mathbf{E}_z[Z_t] = \mathbf{m} + (z - \mathbf{m})e^{t\Phi(1)}.$

(b) $\quad\mathbf{E}_0[Z_t^2] = \dfrac{16(1 - e^{t\Phi(2)})}{\Phi(1)\Phi(2)} + \begin{cases} \dfrac{16t}{\Phi(1)}e^{t\Phi(1)}, & if\ \Phi(1) = \Phi(2),\\[2ex] \dfrac{16(e^{t\Phi(2)} - e^{t\Phi(1)})}{\Phi(1)(\Phi(2) - \Phi(1))}, & otherwise. \end{cases}$

PROOF. $U$ under $\mathbf{P}_z$ is a squared Bessel process of dimension 2 starting from $z$, therefore, $\mathbf{E}_z[U(x)] = z + 2x$. Making use of the independence of $U$ and $\mathbb{V}$, we get

$$\mathbf{E}_z[Z_t] = \mathbf{E}_z[e^{\mathbb{V}_t}U(a(t))] = \mathbf{E}[e^{\mathbb{V}_t}\mathbf{E}_z[U(a(t))|\mathbb{V}]] = \mathbf{E}[e^{\mathbb{V}_t}(z + 2a(t))].$$

We already noticed that time reversal of the Lévy process $\mathbb{V}$ implies that $e^{\mathbb{V}_t} a(t)$ and $\int_0^t e^{\mathbb{V}_s} \, ds$ have the same law, therefore,

$$\mathbf{E}_z[Z_t] = z\mathbf{E}[e^{\mathbb{V}_t}] + 2\int_0^t \mathbf{E}[e^{\mathbb{V}_s}] \, ds = ze^{t\Phi(1)} + \frac{2}{\Phi(1)}(e^{t\Phi(1)} - 1)$$

$$= \mathbf{m} + (z - \mathbf{m})e^{t\Phi(1)}.$$

We now prove (b). First, the scaling property of $U$ shows that, under $\mathbf{P}_0$, the random variables $Z_t$ and $e^{\mathbb{V}_t} a(t) U(1)$ have the same law. Second, $e^{\mathbb{V}_t} a(t)$ and $\int_0^t e^{\mathbb{V}_s} \, ds$ also have the same law. Therefore,

$$(6.1) \qquad \mathbf{E}_0[Z_t^2] = \mathbf{E}_0[U(1)^2]\mathbf{E}\left[\left(\int_0^t e^{\mathbb{V}_s} \, ds\right)^2\right] = 8\mathbf{E}\left[\left(\int_0^t e^{\mathbb{V}_s} \, ds\right)^2\right],$$

where we used the fact that $\mathbf{E}_0[U(1)^2] = 8$ because $U(1)$ under $\mathbf{P}_0$ has an exponential law with mean 2. From a change of variable and making use of the stationarity and the independence of the increments of $\mathbb{V}$, we get

$$\mathbf{E}\left[\left(\int_0^t e^{\mathbb{V}_s} \, ds\right)^2\right] = 2\mathbf{E}\left[\int_0^t e^{\mathbb{V}_x} \int_x^t e^{\mathbb{V}_y} \, dy \, dx\right]$$

$$= 2\int_0^t \mathbf{E}\left[e^{2\mathbb{V}_x} \int_x^t e^{\mathbb{V}_y - \mathbb{V}_x} \, dy\right] dx$$

$$= 2\int_0^t \mathbf{E}[e^{2\mathbb{V}_x}] \int_0^{t-x} \mathbf{E}[e^{\mathbb{V}_y}] \, dy \, dx$$

$$= 2\int_0^t e^{x\Phi(2)} \int_0^{t-x} e^{y\Phi(1)} \, dy \, dx$$

$$= \frac{2(1 - e^{t\Phi(2)})}{\Phi(1)\Phi(2)}$$

$$+ \begin{cases} \dfrac{2t}{\Phi(1)} e^{t\Phi(1)}, & \text{if } \Phi(1) = \Phi(2), \\ \dfrac{2(e^{t\Phi(2)} - e^{t\Phi(1)})}{\Phi(1)(\Phi(2) - \Phi(1))}, & \text{otherwise.} \end{cases}$$

This equality combined with (6.1) completes the proof of (b). $\quad\square$

LEMMA 6.2. *When $\kappa > 2$,*

$$\lim_{\lambda \to 0+} \lambda \mathbf{E}_0\left[\left(\int_0^\infty (Z_t - \mathbf{m})e^{-\lambda t} \, dt\right)^2\right] = \frac{4(\Phi(2) - 4\Phi(1))}{\Phi(1)^3 \Phi(2)}.$$

*This limit is strictly positive because $\Phi$ is a strictly convex function with $\Phi(0) = \Phi(\kappa) = 0$.*

PROOF. We write, for $\lambda > 0$,

$$
\begin{aligned}
\mathbf{E}_0 & \left[ \left( \int_0^\infty (Z_t - \mathbf{m}) e^{-\lambda t}\, dt \right)^2 \right] \\
&= \mathbf{E}_0 \left[ \left( \int_0^\infty Z_t e^{-\lambda t}\, dt \right)^2 \right] - \frac{2\mathbf{m}}{\lambda} \mathbf{E}_0 \left[ \int_0^\infty Z_t e^{-\lambda t}\, dt \right] + \frac{\mathbf{m}^2}{\lambda^2}.
\end{aligned}
\tag{6.2}
$$

Making use of (a) of Lemma 6.1, we find, for any $z \geq 0$,

$$
\begin{aligned}
\mathbf{E}_z \left[ \int_0^\infty Z_t e^{-\lambda t}\, dt \right]
&= \int_0^\infty \mathbf{E}_z[Z_t] e^{-\lambda t}\, dt \\
&= \int_0^\infty (\mathbf{m} + (z - \mathbf{m}) e^{t \Phi(1)}) e^{-\lambda t}\, dt = \frac{z\lambda + 2}{\lambda(\lambda - \Phi(1))}.
\end{aligned}
\tag{6.3}
$$

This equality for $z = 0$ combined with (6.2) yields

$$
\begin{aligned}
\mathbf{E}_0 & \left[ \left( \int_0^\infty (Z_t - \mathbf{m}) e^{-\lambda t}\, dt \right)^2 \right] \\
&= \mathbf{E}_0 \left[ \left( \int_0^\infty Z_t e^{-\lambda t}\, dt \right)^2 \right] + \frac{4(\lambda + \Phi(1))}{\lambda^2 \Phi(1)^2 (\lambda - \Phi(1))}.
\end{aligned}
\tag{6.4}
$$

We also have

$$
\begin{aligned}
\mathbf{E}_0 \left[ \left( \int_0^\infty Z_t e^{-\lambda t}\, dt \right)^2 \right]
&= 2 \mathbf{E}_0 \left[ \int_0^\infty Z_x e^{-\lambda x} \int_x^\infty Z_y e^{-\lambda y}\, dy\, dx \right] \\
&= 2 \int_0^\infty \mathbf{E}_0 \left[ Z_x e^{-\lambda x} \int_0^\infty Z_{x+y} e^{-\lambda(x+y)}\, dy \right] dx \\
&= 2 \int_0^\infty e^{-2\lambda x} \mathbf{E}_0 \left[ Z_x \mathbb{E}_{Z_x} \left[ \int_0^\infty Z_y e^{-\lambda y}\, dy \right] \right] dx,
\end{aligned}
$$

where we used the Markov property of $Z$ for the last equality. Thus, with the help of (6.3), we find

$$
\mathbf{E}_0 \left[ \left( \int_0^\infty Z_t e^{-\lambda t}\, dt \right)^2 \right] = \frac{2}{\lambda(\lambda - \Phi(1))} \int_0^\infty e^{-2\lambda x} (\lambda \mathbf{E}_0[Z_x^2] + 2 \mathbf{E}_0[Z_x])\, dx.
$$

This integral can now be explicitly computed thanks to Lemma 6.1. After a few lines of elementary calculus, we obtain

$$
\mathbf{E}_0 \left[ \left( \int_0^\infty Z_t e^{-\lambda t}\, dt \right)^2 \right] = \frac{4(6\lambda - \Phi(2))}{\lambda^2 (\lambda - \Phi(1))(4\lambda^2 - 2\lambda(\Phi(1) + \Phi(2)) + \Phi(1)\Phi(2))}
$$

[this result does not depend on whether or not $\Phi(1) = \Phi(2)$]. Substituting this equality in (6.4), we get

$$
\lambda \mathbf{E}_0 \left[ \left( \int_0^\infty (Z_t - \mathbf{m}) e^{-\lambda t}\, dt \right)^2 \right]
$$

$$= \frac{4\Phi(1)(\Phi(2) - 4\Phi(1)) - 4\lambda(4\lambda + 2(\Phi(1) - \Phi(2)))}{\Phi(1)^2(\Phi(1) - \lambda)(4\lambda^2 - 2\lambda(\Phi(1) + \Phi(2)) + \Phi(1)\Phi(2))}.$$

We conclude the proof of the lemma by taking the limit as $\lambda$ tends to $0+$.
□

LEMMA 6.3.    *When* $\kappa > 2$,

$$\lim_{\lambda \to 0+} \lambda \mathbf{E}_1 \left[ \left( \int_0^\infty (Z_t - \mathbf{m}) e^{-\lambda t} \, dt \right)^2 \right] = \frac{1}{2\mathbf{n}[\zeta]} \mathbf{n} \left[ \left( \int_0^\zeta (\epsilon_t - \mathbf{m}) \, dt \right)^2 \right].$$

PROOF.    Recall that $\varphi$ stands for the Laplace exponent of the inverse of the local time $L^{-1}$. We first use (b) of Lemma 5.2 with the function $f(x) = |x - \mathbf{m}|$:

$$\mathbf{E}_1 \left[ \left( \int_0^\infty |Z_t - \mathbf{m}| e^{-\lambda t} \, dt \right)^2 \right]$$

$$(6.5) \qquad = \frac{1}{\varphi(2\lambda)} \mathbf{n} \left[ \left( \int_0^\zeta |\epsilon_t - \mathbf{m}| e^{-\lambda t} \, dt \right)^2 \right]$$

$$\qquad\qquad + \frac{2}{\varphi(\lambda)\varphi(2\lambda)} \mathbf{n} \left[ \int_0^\zeta |\epsilon_t - \mathbf{m}| e^{-\lambda t} \, dt \right] \mathbf{n} \left[ e^{-\lambda \zeta} \int_0^\zeta |\epsilon_t - \mathbf{m}| e^{-\lambda t} \, dt \right].$$

Note also that Lemma 5.1 and Proposition 5.2 readily show that

$$(6.6) \qquad \mathbf{n} \left[ \left( \int_0^\zeta |\epsilon_t - \mathbf{m}| \, dt \right)^\beta \right] < \infty \qquad \text{for all } \beta < \kappa.$$

Thus, the three expectations under $\mathbf{n}$ on the r.h.s. of (6.5) are finite because $\kappa > 2$. Therefore, we can also use (b) of Lemma 5.2 with the function $f(x) = x - \mathbf{m}$:

$$\mathbf{E}_1 \left[ \left( \int_0^\infty (Z_t - \mathbf{m}) e^{-\lambda t} \, dt \right)^2 \right]$$

$$(6.7) \qquad = \frac{1}{\varphi(2\lambda)} \mathbf{n} \left[ \left( \int_0^\zeta (\epsilon_t - \mathbf{m}) e^{-\lambda t} \, dt \right)^2 \right]$$

$$\qquad\qquad + \frac{2}{\varphi(\lambda)\varphi(2\lambda)} \mathbf{n} \left[ \int_0^\zeta (\epsilon_t - \mathbf{m}) e^{-\lambda t} \, dt \right] \mathbf{n} \left[ e^{-\lambda \zeta} \int_0^\zeta (\epsilon_t - \mathbf{m}) e^{-\lambda t} \, dt \right].$$

Recall also that $\varphi(\lambda) \sim \mathbf{n}[\zeta]\lambda$ [cf. (5.7) in the proof of Corollary 5.1]. Thus, keeping in mind (6.6), the dominated convergence theorem yields

$$\lim_{\lambda \to 0+} \frac{\lambda}{\varphi(2\lambda)} \mathbf{n} \left[ \left( \int_0^\zeta (\epsilon_t - \mathbf{m}) e^{-\lambda t} \, dt \right)^2 \right] = \frac{1}{2\mathbf{n}[\zeta]} \mathbf{n} \left[ \left( \int_0^\zeta (\epsilon_t - \mathbf{m}) \, dt \right)^2 \right]$$

and

$$\lim_{\lambda \to 0+} \mathbf{n}\left[ e^{-\lambda \zeta} \int_0^\zeta (\epsilon_t - \mathbf{m}) e^{-\lambda t} \, dt \right] = \mathbf{n}\left[ \int_0^\zeta (\epsilon_t - \mathbf{m}) \, dt \right] = 0,$$

where we used Corollary 5.2 for the last equality. Finally, (a) of Lemma 5.2 combined with (6.3) gives

$$\frac{2\lambda}{\varphi(\lambda)\varphi(2\lambda)} \mathbf{n}\left[ \int_0^\zeta (\epsilon_t - \mathbf{m}) e^{-\lambda t} \, dt \right]$$

$$= \frac{2\lambda}{\varphi(2\lambda)} \mathbf{E}_1\left[ \int_0^\infty (Z_t - \mathbf{m}) e^{-\lambda t} \, dt \right]$$

$$= \frac{2\lambda}{\varphi(2\lambda)} \left( \frac{\Phi(1) + 2}{\Phi(1)(\lambda - \Phi(1))} \right) \underset{\lambda \to 0+}{\longrightarrow} - \frac{2 + \Phi(1)}{\mathbf{n}[\zeta]\Phi(1)^2}.$$

These last three estimates combined with (6.7) entail the lemma.   $\square$

We can now easily obtain the calculation of the second moment.

PROPOSITION 6.1.   *When $\kappa > 2$,*

$$\mathbf{n}\left[ \left( \int_0^\zeta (\epsilon_t - \mathbf{m}) \, dt \right)^2 \right] = \mathbf{n}[\zeta] \frac{8(\Phi(2) - 4\Phi(1))}{\Phi(1)^3 \Phi(2)}.$$

PROOF.   In view of Lemmas 6.2 and 6.3, it suffices to prove that

$$\lim_{\lambda \to 0+} \lambda \mathbf{E}_1\left[ \left( \int_0^\infty (Z_t - \mathbf{m}) e^{-\lambda t} \, dt \right)^2 \right] = \lim_{\lambda \to 0+} \lambda \mathbf{E}_0\left[ \left( \int_0^\infty (Z_t - \mathbf{m}) e^{-\lambda t} \, dt \right)^2 \right].$$

Indeed, the Markov property of $Z$ for the stopping time $\tau_1$ yields

$$\mathbf{E}_0\left[ \left( \int_0^\infty (Z_t - \mathbf{m}) e^{-\lambda t} \, dt \right)^2 \right]$$

$$= \mathbf{E}_0\left[ \left( \int_0^{\tau_1} (Z_t - \mathbf{m}) e^{-\lambda t} \, dt + \int_{\tau_1}^\infty (Z_t - \mathbf{m}) e^{-\lambda t} \, dt \right)^2 \right]$$

(6.8)   $$= \mathbf{E}_0\left[ \left( \int_0^{\tau_1} (Z_t - \mathbf{m}) e^{-\lambda t} \, dt \right)^2 \right]$$

$$\quad + \mathbf{E}_0[e^{-2\lambda \tau_1}] \mathbf{E}_1\left[ \left( \int_0^\infty (Z_t - \mathbf{m}) e^{-\lambda t} \, dt \right)^2 \right]$$

$$\quad + 2\mathbf{E}_0\left[ e^{-\lambda \tau_1} \int_0^{\tau_1} (Z_t - \mathbf{m}) e^{-\lambda t} \, dt \right] \mathbf{E}_1\left[ \int_0^\infty (Z_t - \mathbf{m}) e^{-\lambda t} \, dt \right].$$

Proposition 4.1 and the absence of positive jumps for $Z$ give

(6.9)   $$\mathbf{E}_0\left[ \left( \int_0^{\tau_1} (Z_t - \mathbf{m}) e^{-\lambda t} \, dt \right)^2 \right] \leq (\mathbf{m} + 1)^2 \mathbf{E}_0[\tau_1^2] < \infty.$$

Similarly,

$$(6.10) \qquad \left| \mathbf{E}_0 \left[ e^{-\lambda \tau_1} \int_0^{\tau_1} (Z_t - \mathbf{m}) e^{-\lambda t} \, dt \right] \right| \leq (\mathbf{m} + 1) \mathbb{E}_0[\tau_1] < \infty.$$

Note also that, according to (6.3),

$$(6.11) \qquad \lambda \mathbf{E}_1 \left[ \int_0^\infty (Z_t - \mathbf{m}) e^{-\lambda t} \, dt \right] = \frac{\lambda(\Phi(1) + 2)}{\Phi(1)(\lambda - \Phi(1))} \xrightarrow[\lambda \to 0+]{} 0.$$

Thus, (6.8)–(6.11) and the fact that $\lim_{\lambda \to 0+} \mathbf{E}_0[e^{-2\lambda \tau_1}] = 1$ conclude the proof of the proposition. $\square$

## 7. End of the proof of the main theorem.

We showed in Section 2 that we simply need to prove Theorem 1.1 for the additive functional

$$I(r) \stackrel{\text{def}}{=} \int_0^r Z_s \, ds$$

under $\mathbf{P} = \mathbf{P}_0$ in place of $H(r)$. Moreover, Proposition 4.3 states that the hitting time of level 1 by $Z$ is $\mathbf{P}_0$-almost surely finite, therefore, it is sufficient to prove the result for $I(r)$ under $\mathbf{P}_1$. The remaining portion of the proof is now quite standard and very similar to the argument given on pages 166–167 of [14]. Let us first deal with the case $\kappa < 1$. Recall that $L^{-1}$ stands for the inverse of the local time of $Z$ at level 1. Since $I$ is an additive functional of $Z$, the process $(I(L_t^{-1}), t \geq 0)$ under $\mathbf{P}_1$ is a subordinator (without drift thanks to Lemma 3.3) whose Laplace transform is given by

$$(7.1) \qquad \mathbf{E}_1[e^{-\lambda I(L_t^{-1})}] = \exp\left( -t\lambda \int_0^\infty e^{-\lambda x} \mathbf{n}\{\widetilde{I}(\epsilon) > x\} \, dx \right),$$

where we used the notation $\widetilde{I}(\epsilon) \stackrel{\text{def}}{=} \int_0^\zeta \epsilon_t \, dt$. We now define

$$\xi_n \stackrel{\text{def}}{=} \int_{L_{n-1}^{-1}}^{L_n^{-1}} Z_s \, ds = I(L_n^{-1}) - I(L_{n-1}^{-1}).$$

The sequence $(\xi_n, n \geq 1)$ under $\mathbf{P}_1$ is i.i.d. Moreover, in view of Proposition 5.2, we deduce from (7.1) that

$$(7.2) \qquad \mathbf{P}_1\{\xi_1 > x\} \underset{x \to \infty}{\sim} \mathbf{n}\{\widetilde{I}(\epsilon) > x\} \underset{x \to \infty}{\sim} \mathbf{n}[\zeta] \frac{2^\kappa \Gamma(\kappa) \kappa^2 \mathbf{K}^2}{\Phi'(\kappa) x^\kappa}.$$

The characterization of the domains of attraction to a stable law (see, e.g., Chapter IX.8 of [9]) implies that

$$\frac{I(L_n^{-1})}{n^{1/\kappa}} = \frac{\xi_1 + \cdots + \xi_n}{n^{1/\kappa}} \xrightarrow[n \to \infty]{\text{law}} 2 \left( \frac{\mathbf{n}[\zeta] \pi \kappa^2 \mathbf{K}^2}{2 \sin(\pi \kappa/2) \Phi'(\kappa)} \right)^{1/\kappa} \mathcal{S}_\kappa^{\text{ca}}.$$

Moreover, according to Lemma 5.1, we have $\mathbf{E}[L_1^{-1}] = \mathbf{n}[\zeta] < \infty$ so the strong law of large numbers for subordinators (cf. page 92 of [1]) yields

$$(7.3) \qquad \frac{L_t^{-1}}{t} \xrightarrow[t\to\infty]{\text{a.s.}} \mathbf{n}[\zeta].$$

We can therefore use Theorem 8.1 of [19] with the change of time $L^{-1}$ to check that, under $\mathbf{P}_1$,

$$\frac{I(t)}{t^{1/\kappa}} \xrightarrow[t\to\infty]{\text{law}} 2\left(\frac{\pi\kappa^2\mathbf{K}^2}{2\sin(\pi\kappa/2)\Phi'(\kappa)}\right)^{1/\kappa} \mathcal{S}_\kappa^{\text{ca}}.$$

This concludes the proof of the theorem when $\kappa < 1$. Let us now assume that $\kappa = 1$. In this case, $\mathbf{K} = \mathbf{E}[(\int_0^\infty e^{\mathbb{V}_s}\,ds)^0] = 1$, hence, (7.2) takes the form

$$\mathbf{P}_1\{\xi_1 > x\} \underset{x\to\infty}{\sim} \frac{2\mathbf{n}[\zeta]}{\Phi'(1)x}.$$

The characterization of the domains of attraction now states that there exists a constant $c_{28}$ such that

$$(7.4) \qquad \frac{I(L_n^{-1}) - ng(n)}{n} = \frac{\xi_1 + \cdots + \xi_n}{n} - g(n) \xrightarrow[n\to\infty]{\text{law}} c_{28} + \frac{\pi\mathbf{n}[\zeta]}{\Phi'(1)}\mathcal{C}^{\text{ca}},$$

where $g(x) \overset{\text{def}}{=} \int_0^x \mathbf{P}_1\{\xi_1 > y\}\,dy$. Note also that our estimate on $\mathbf{n}\{\zeta > x\}$ (Lemma 5.1) entails an iterated logarithm law for the subordinator $L^{-1}$, in particular,

$$\frac{L_n^{-1}}{\mathbf{n}[\zeta]} \in [n - n^{2/3}, n + n^{2/3}] \qquad \text{for all } n \text{ large enough.}$$

Using this result and the fact that $I(\cdot)$ is nondecreasing, it is not difficult to deduce from (7.4) that

$$\frac{1}{t}\left(I(t) - \frac{t}{\mathbf{n}[\zeta]}g\left(\frac{t}{\mathbf{n}[\zeta]}\right)\right) \xrightarrow[t\to\infty]{\text{law}} \frac{c_{28}}{\mathbf{n}[\zeta]} + \frac{\pi}{\Phi'(1)}\mathcal{C}^{\text{ca}}$$

(cf. with the argument given on page 166 of [14] for details). Thus, setting

$$f(t) \overset{\text{def}}{=} \frac{t}{\mathbf{n}[\zeta]}\left(g\left(\frac{t}{\mathbf{n}[\zeta]}\right) - c_{28}\right),$$

we get the desired limiting law for $(I(t) - f(t))/t$ and also

$$f(t) \underset{t\to\infty}{\sim} \frac{t}{\mathbf{n}[\zeta]}\int_0^{t/(\mathbf{n}[\zeta])} \mathbf{P}_1\{\xi_1 > y\}\,dy \underset{t\to\infty}{\sim} \frac{2t\log t}{\Phi'(1)}.$$

The proof of the theorem when $\kappa > 1$ is very similar to that in the case $\kappa < 1$, but we now consider the sequence $(\xi_n', n \geq 1)$ instead of $(\xi_n, n \geq 1)$ defined by

$$\xi_n' \overset{\text{def}}{=} \int_{L_{n-1}^{-1}}^{L_n^{-1}} (Z_s - \mathbf{m})\,ds = \xi_n - \mathbf{m}(L_n^{-1} - L_{n-1}^{-1}).$$

These random variables are i.i.d. and are centered under $\mathbf{P}_1$ because

$$\mathbf{E}_1[\xi_1'] = \mathbf{n}\left[\int_0^\zeta (\epsilon_s - \mathbf{m})\,ds\right] = \mathbf{n}\left[\int_0^\zeta \epsilon_s\,ds\right] - \mathbf{n}[\zeta]\mathbf{m} = 0$$

(we used Corollary 5.2 for the last equality). Moreover, when $\kappa > 2$, Proposition 6.1 yields

$$\mathbf{E}_1[{\xi_1'}^2] = \mathbf{n}\left[\left(\int_0^\zeta (\epsilon_s - \mathbf{m})\,ds\right)^2\right] = \mathbf{n}[\zeta]\frac{8(\Phi(2) - 4\Phi(1))}{\Phi(1)^3\Phi(2)}.$$

Since the tail distribution of $\zeta$ under $\mathbf{n}$ has (at least) an exponential decrease, we see that the estimate (7.2) still holds with $\xi_1'$ in place of $\xi_1$. Thus, the characterization of the domains of attraction to a stable law insures that, when $\kappa \in (1, 2)$,

$$\frac{I(L_n^{-1}) - \mathbf{m}L_n^{-1}}{n^{1/\kappa}} = \frac{\xi_1' + \cdots + \xi_n'}{n^{1/\kappa}}\, n \to \infty \xrightarrow{\text{law}} 2\left(\frac{\mathbf{n}[\zeta]\pi\kappa^2\mathbf{K}^2}{2\sin((\pi\kappa)/2)\Phi'(\kappa)}\right)^{1/\kappa}\mathcal{S}_\kappa^{\text{ca}}.$$

Similarly, when $\kappa = 2$ and since $\mathbf{K} = \mathbf{E}[\int_0^\infty e^{\mathbb{V}_s}\,ds] = \frac{-1}{\Phi(1)}$,

$$\frac{I(L_n^{-1}) - \mathbf{m}L_n^{-1}}{\sqrt{n\log n}} \xrightarrow[n\to\infty]{\text{law}} \frac{-4\sqrt{\mathbf{n}[\zeta]}}{\Phi(1)\sqrt{\Phi'(2)}}\mathcal{N},$$

and when $\kappa > 2$,

$$\frac{I(L_n^{-1}) - \mathbf{m}L_n^{-1}}{\sqrt{n}} \xrightarrow[n\to\infty]{\text{law}} \sqrt{\mathbf{E}_1[{\xi_1'}^2]}\mathcal{N} = \sqrt{\frac{\mathbf{n}[\zeta]8(\Phi(2) - 4\Phi(1))}{\Phi(1)^3\Phi(2)}}\mathcal{N}.$$

Just as in the case $\kappa < 1$, we easily check that the hypotheses of Theorem 8.1 of [19] are fulfilled. Thus, the change of time $L_t^{-1} \sim \mathbf{n}[\zeta]t$ is legitimate and concludes the proof of the theorem. $\square$

**Acknowledgment.** I am most grateful to my Ph.D. supervisor Yueyun Hu for his unfailing help and advice.

LABORATOIRE DE PROBABILITÉS
ET MODÈLES ALÉATOIRES
UNIVERSITÉ PIERRE ET MARIE CURIE
175 RUE DU CHEVALERET
75013 PARIS
FRANCE
E-MAIL: arvind.singh@ens.fr